\documentclass{article}

\usepackage{amssymb}
\usepackage{amsfonts}
\usepackage{amsmath}
\usepackage{url}
\usepackage{latexsym}

\def\deg{\mathop{\rm deg }\nolimits}
\def\rank{\mathop{\rm rank}\nolimits}
\def\lcm{\mathop{\rm lcm }\nolimits}
\def\rev{\mathop{\rm rev }\nolimits}
\newcommand{\se}{\ensuremath{\stackrel{s.e.}{\sim}}}

\def\gr{\mathop{\rm grade }\nolimits}

\newcommand{\efe}{\mathbb F}
\newcommand{\FF}{\mathbb F}

\newcommand{\la}{s}

\newcommand{\ba}{\mathbf a}
\newcommand{\bd}{\mathbf d}

\newcommand{\bc}{\mathbf c}
\newcommand{\bb}{\mathbf b}
\newcommand{\bu}{\mathbf u}
\newcommand{\bv}{\mathbf v}

\newtheorem{theorem}{Theorem}[section]
\newtheorem{corollary}[theorem]{Corollary}
\newtheorem{lemma}[theorem]{Lemma}
\newtheorem{definition}[theorem]{Definition}
\newtheorem{problem}[theorem]{Problem}
\newtheorem{proposition}[theorem]{Proposition}
\newtheorem{remark}[theorem]{Remark}
\newtheorem{example}[theorem]{Example}
\newenvironment{rem}{\begin{remark} \em}{\end{remark}}

\title{Row completion of polynomial and rational matrices
  \thanks{
This work was 
supported by grants PID2021-124827NB-I00 and RED2022-134176-T funded by MCIN/AEI/ 10.13039/501100011033 and by ``ERDF A way of making Europe'' by the ``European Union''.
The first and third authors were also supported 
by grant GIU21/020 funded by UPV/EHU.
  }
  }

\author{Agurtzane Amparan\thanks{Departamento de Matem\'aticas, Universidad del Pa\'is Vasco UPV/EHU, Bilbao, Spain, {agurtzane.amparan@ehu.eus}, {silvia.marcaida@ehu.eus}}
\and Itziar Baraga\~na\thanks{Departamento de Ciencia de la Computaci\'on e I.A., 
Universidad del Pa\'{\i}s Vasco, UPV/EHU, Donostia-San Sebasti\'an, 
Spain, {itziar.baragana@ehu.eus}}
\and Silvia Marcaida\footnotemark[2]
\and Alicia Roca\thanks{Departamento de Matem\'atica Aplicada, IMM, Universitat Polit\`ecnica de Val\`encia, 46022 Valencia, Spain,   {aroca@mat.upv.es}}}

\date{}

\DeclareMathOperator{\diag}{diag}

\begin{document}

\maketitle

\begin{abstract}
  We characterize the existence of a polynomial (rational) matrix when its eigenstructure (complete structural data) and 
some of its rows are prescribed. For polynomial matrices, this problem was solved in \cite{AAmBaMaRo24} when the polynomial matrix has the same degree as the prescribed submatrix. In that paper,  the following row completion problems  were also solved arising  when the eigenstructure was partially prescribed, keeping the  restriction on the degree: the eigenstructure but the row
 (column) minimal indices, and the finite and/or infinite structures.
 Here we remove the restriction on the degree, allowing it to be greater than or equal to that of the submatrix. We also generalize the results to rational matrices. Obviously, the results obtained hold for the corresponding column  completion problems.
\end{abstract}


{\bf AMS:}
 15A54, 15A83, 93B18

\section{Introduction}

An important problem in Matrix Theory is the matrix completion problem. It  consists in characterizing the existence of a matrix with certain properties when a submatrix is prescribed. In fact,  this  problem includes many other ones depending on the type  of  matrices involved and the properties analyzed. In the last decades the research in the area has been very fruitful. See \cite{AAmBaMaRo24, DoSt19, LoMoZaZaLAA98, Sa79, Th79}  and the references therein.

This work is devoted to the  matrix completion problem for polynomial and rational matrices when  the   complete structural data (or some of them) of the polynomial or the rational matrix are prescribed and the submatrix is formed by some of its rows (columns). This study generalizes  the results obtained  in \cite{AAmBaMaRo24}, where the row completion problem of a polynomial matrix is solved when the eigenstructure (or part of it) is prescribed and the degree of the completed matrix is the same as that of the prescribed submatrix. 

The generalization addressed in this paper is two-folded. On the one hand, we allow that the degree of the completed polynomial matrix is greater than or equal to that of the prescribed submatrix. On the other hand, the results of \cite{AAmBaMaRo24} are generalized to rational matrices, solving  the row completion problem when  the complete structural data (or some of them) of the completed rational matrix are prescribed. 

The eigenstructure of a polynomial matrix is formed by four types of invariants: the invariant factors, the partial multiplicities of $\infty$, and the column and row minimal indices (\cite{VDoDe83}). The invariant factors form the finite structure of the matrix, the partial multiplicities of $\infty$ are the infinite structure, and the minimal indices, the singular structure.  The prescription of all or part of these invariants leads to pose 15 different completion problems, some of them solved in \cite{AAmBaMaRo24}, when the degree of the completed matrix is prescribed to  coincide with that of the given  submatrix.  The remaining ones are solved in \cite{2AmBaMaRo24}.

The complete structural data of a rational matrix is formed by the invariant rational functions (finite structure), the invariant orders at $\infty$ (infinite structure), and the column and row minimal indices (singular structure) (\cite{AnDoHoMa19}). 
We show later that for polynomial matrices, knowing the eigenstructure is equivalent to knowing the complete structural data.

Once more, due to the high number of problems that the complete analysis of the study includes,  we solve in this work some of the cases, generalizing the results in \cite{AAmBaMaRo24}, and leaving the remaining ones for a future paper, hence generalizing \cite{2AmBaMaRo24}.
One of the problems, the  prescription of only the finite structure, was solved independently by E. Marques de S\'a  and R. C. Thompson in 1979 in the seminal papers \cite{Sa79,Th79}, and the result was generalized to the rational case  in \cite{Ba89}.

The paper is organized as follows. 
Section \ref{secprelimin} contains the notation, definitions, preliminary results,
and the statement of the problems we deal with. In Section \ref{secpolcompl} we solve the problem when the complete structural data are prescribed. Finally, in Section \ref{secpart} we solve the problem of the prescription of the complete structural data  but the row (column) minimal indices (Subsection \ref{subsec_fininfcolrow}), 
and that of the prescription of the finite and/or infinite structures (Subsection \ref{subsec_dinoronf}). 

\section{Preliminaries}\label{secprelimin}

Let $\FF$ be a field. 
The ring of polynomials in the indeterminate $s$ with coefficients in $\FF$
is denoted by $\FF[s]$,  $\FF(s)$ is the field of fractions of
$\FF[s]$, i.e., the field of rational functions over $\FF$, and $\FF_{pr}(s)$ is the ring of proper rational functions, i.e., the rational functions with degree of the denominator at least the degree of the numerator. The ring of polynomials in two
variables $s, t$ with coefficients in $\FF$
is denoted by $\FF[s, t]$. 
A polynomial in $\FF[s]$ is \textit{monic} if its leading coefficient is 1. We say that a polynomial in $\FF[s,t]$ is \textit{monic}
    if it is monic with respect to the variable $s$.
  Given two polynomials $\alpha, \beta$, by $\alpha\mid \beta$ we mean that $\alpha$ is a divisor of $\beta$, by $\lcm(\alpha,\beta)$, the monic least common multiple of $\alpha$ and $\beta$, and by $\gcd(\alpha,\beta)$, the monic greatest common divisor of $\alpha$ and $\beta$.

   In this work we deal with  finite sequences of integers $\ba=(a_1, \dots, a_r)$ where   $a_1\geq\dots\geq a_r$. If $a_r\geq 0$,
the sequence is called a 
\textit{partition}. When necessary, we take $a_i=+\infty$ for $i<1$ and 
$a_i=-\infty$ for $i>r$.  If $b_1\leq\dots\leq b_r$ is an increasing sequence of integers,  we take $b_i=-\infty$ for $i<1$ and 
$b_i=+\infty$ for $i>r$. 

 We also deal with  polynomial chains $\alpha_1\mid \dots \mid \alpha_r$, where $\alpha_i\in \FF[s]$ or
 $\alpha_i\in \FF[s,t]$, and take $\alpha_1=1$ for $i<1$ and $\alpha_i=0$ for $i>r$.
 If $\varphi_r\mid \dots \mid \varphi_1$,  we take $\varphi_i=1$ for $i>r$ and $\varphi_i=0$ for $i<1$.

We denote by  $\FF^{m\times n}$, $\FF[s]^{m\times n}$, $\FF(s)^{m \times n}$, and $\FF_{pr}(s)^{m \times n}$ the vector spaces over $\FF$  of $m \times n$ matrices with elements in $\FF$, $\FF[s]$, $\FF(s)$, and $\FF_{pr}(s)$, respectively.
A matrix $U(s)\in\efe[s]^{n\times n}$ is said \textit{unimodular} if it has inverse in $\efe[s]^{n\times n}$, while a matrix $B(s)\in\efe_{pr}(s)^{n\times n}$ is said \textit{biproper} if it has inverse in $\efe_{pr}(s)^{n\times n}$.

Let $R(s)\in \FF(s)^{m\times n}$ of  $\rank(R(s))=r$. A canonical form for the unimodular equivalence of $R(s)$   is the \textit{Smith--McMillan form}
$$
\begin{bmatrix}\diag\left(\frac{\eta_1(s)}{\varphi_1(s)},\dots,\frac{\eta_r(s)}{\varphi_r(s)}\right)&0\\0&0\end{bmatrix},$$
where $\eta_1(s)\mid  \dots \mid \eta_r(s)$ and $\varphi_r(s)\mid  \dots \mid \varphi_1(s)$ are monic polynomials and $\frac{\eta_1(s)}{\varphi_1(s)},\dots,\frac{\eta_r(s)}{\varphi_r(s)}$ are irreducible rational functions
 called the \textit{invariant rational functions} of $R(s)$. We also refer to them as the \textit{finite structure} of $R(s)$. The polynomial $\varphi_1(s)$ is the monic least common denominator of the entries of $R(s)$ (see, for instance, \cite[Chapter 3, Section 4]{Rose70}).

A canonical form for the equivalence at infinity of $R(s)$ is the  \textit{Smith--McMillan  form at infinity}  
$$ 
\begin{bmatrix}\diag\left(s^{-\tilde p_1},\ldots,s^{-\tilde p_r}\right)&0\\0&0\end{bmatrix},
$$
where  $\tilde p_1\leq\cdots\leq\tilde p_r$ are integers called the \textit{invariant orders at infinity} of $R(s)$  (see, for instance, \cite{Vard91}).
In \cite{AnDoHoMa19}, the sequence  of invariant orders at $\infty$  is called  the \textit{structural index sequence of $R(s)$ at $\infty$.}

We recall now the \textit{singular structure} of a rational matrix. Denote by $\mathcal{N}_\ell (R(\la))$ and $\mathcal{N}_r (R(\la))$ the \textit{left} and
\textit{right null-spaces} over $\FF(\la)$ of $R(\la)$, respectively, i.e.,
if $R(\la)\in\FF(\la)^{m\times n}$,
\[
\begin{array}{l}
\mathcal{N}_\ell (R(\la))=\{x(\la)\in\FF(\la)^{m \times 1}: x(\la)^TR(\la)=0\},\\
\mathcal{N}_r (R(\la))=\{x(\la)\in\FF(\la)^{n \times 1}: R(\la)x(\la)=0\},
\end{array}
\]
which are vector subspaces
of $\FF(\la)^{m \times 1}$ and $\FF(\la)^{n \times 1}$, respectively. For a subspace $\mathcal{V}$ of
$\FF(\la)^{m \times 1}$ it is possible to find a basis consisting of vector
polynomials; it is enough to  take an arbitrary basis and multiply each vector
by a  least common multiple of the denominators of its  entries.
The \textit{order} of a polynomial basis is defined as the sum of
the degrees of its vectors  (see \cite{Fo75}). A \textit{minimal basis} of  $\mathcal{V}$
is  a polynomial basis  with least order among the
polynomial bases of $\mathcal{V}$. 
The  degrees of the vector polynomials  of a minimal basis, increasingly ordered,  are  always the same (see \cite{Fo75}), and are called the \textit{minimal indices} of $\mathcal{V}$.

 A   \textit{right  (left) minimal basis} of a rational matrix $R(\la)$ is a minimal basis of $\mathcal{N}_r (R(\la))$ ($\mathcal{N}_\ell (R(\la))$).
The \textit{right  (left) minimal indices} of $R(\la)$ are the minimal indices of
$\mathcal{N}_r (R(\la))$ ($\mathcal{N}_\ell (R(\la))$). 
From now on in this paper, we  work with the right (left) minimal indices decreasingly ordered, and we  refer to them as the \textit{column (row) minimal indices} of $R(s)$.  Notice that a rational matrix $R(s)\in \FF(s)^{m\times n}$ of $\rank(R(s))=r$ has $m-r$ row  and $n-r$ column minimal indices.

 For a rational matrix $R(s)\in \FF(s)^{m\times n}$
of $\rank(R(s))=r$, the \textit{complete structural data} consist of four components
(see \cite[Definition 2.15]{AnDoHoMa19}): the invariant rational functions $\frac{\eta_1(s)}{\varphi_1(s)},\dots,   \frac{\eta_r(s)}{\varphi_r(s)}$,
the  invariant orders at infinity  $\tilde p_1\leq \cdots \leq  \tilde p_{r}$,  the row minimal indices $(u_1, \dots, u_{m-r})$ and the column minimal indices $(c_1, \dots, c_{n-r})$.
Observe that the complete structural data of a rational matrix determine its $\rank$.

 If the rational matrix is a polynomial matrix $P(s)$ of $\rank P(s)=r$,  then $\varphi_1(s)=\dots=\varphi_r(s)=1$, the polynomials $\eta_1(s)\mid\dots\mid\eta_r(s)$ are the \textit{invariant factors} of $R(s)$, and the Smith-McMillan form is its \textit{Smith normal form} (\cite[Chapter 1, Section 1]{Rose70}).
Hence, the complete structural data of a polynomial matrix is formed by the invariant factors, the invariant orders at infinity, and the column and row minimal indices.

For polynomial matrices we introduce some other definitions. 
Let  $\deg (P(s))=d$, where $\deg(\cdot)$ stands for degree.
The \textit{reversal} of $P(s)$ is the polynomial matrix
  $$\rev(P)(t)=t^d P\left(\frac{1}{t}\right).$$
The \textit{partial multiplicities of $\infty$} in $P(s)$ are defined as the partial multiplicities of $0$ in  $\rev(P)(t)$ (see, for instance,   \cite{DeDoVa15}). 

The invariant factors, the  partial multiplicities of $\infty$,  the row minimal indices and the column minimal indices  are known as the \textit{eigenstructure} of the polynomial matrix $P(s)\in \FF[s]^{m\times n}$
(see \cite{VDoDe83}).

Observe that the eigenstructure of a polynomial matrix determines its $\rank$ (it is the number of  invariant factors, or the number of  partial multiplicities of $\infty$, and it is also equal to  the number of columns (rows) minus the number of column (row) minimal indices).

In the literature, the invariant factors and the partial multiplicities of $\infty$ of a polynomial matrix $P(s)$ are often treated together as follows:
  Let  $\alpha_1(s)\mid  \dots \mid \alpha_r(s)$ be the  invariant factors  and $e_1\leq \dots \leq e_r$ the partial multiplicities of $\infty$ of $P(s)$.
The \textit{homogeneous invariant factors  of $P(s)$}
are homogeneous polynomials in $\FF[s,t]$,
$\phi_{1}(s,t) \mid \dots \mid \phi_{r}(s,t)$, defined as 
$$
  \phi_{i}(s,t)=t^{e_{i}}t^{\deg(\alpha_i)}\alpha_i\left(\frac st\right), \quad 1\leq i\leq r.
  $$

Given $P(s), \bar P(s)\in\efe[s]^{m\times n}$, we write $P(s)\approx \bar  P(s)$ when they have the same eigenstruture. If  $P(s), \bar P(s)$ are matrix pencils,  $P(s)\approx \bar  P(s)$ if and only if  they are strictly equivalent ($P(s)\se\bar  P(s)$), i.e., $\bar P(s)=S P(s)T$ for some non singular matrices $S$ and  $T$. 

 Prior to state our problems, we present some results related to the existence of polynomial or rational matrices with prescribed eigenstructure or complete structural data.

\begin{theorem}[{\rm \cite[Theorem 3.1] {AAmBaMaRo24}, \cite[Theorem 3.3]{DeDoVa15} for infinite fields}]\label{theoexistenceDeDoVa152}
Let $m$, $n$, $r$ be positive integers, $r\leq \min\{m,n\}$, and $d$ a non negative integer.
Let
$\alpha_1(s)\mid \dots \mid \alpha_r(s)$ be monic  polynomials.
Let $(e_r, \ldots, e_{1})$, $(c_1, \ldots, c_{n-r})$, $(u_1, \ldots, u_{m-r})$
be partitions. 
Then, there exists a polynomial matrix  $P(s)\in \FF[s]^{m\times n}$ of $\rank(P(s))= r$,
$\deg(P(s))= d$,  with  $\alpha_1(s), \dots, \alpha_r(s)$ as invariant factors, $e_1, \dots, e_r$ as partial multiplicities of $\infty$, and $c_1, \dots, c_{n-r}$ and $u_1, \dots, u_{m-r}$ as column and row minimal indices,
 respectively, if and only if
 \begin{equation*}\label{eqf1}
  e_1= 0,
\end{equation*}
\begin{equation}\label{eqIST}
\sum_{i=1}^{n-r}c_i+\sum_{i=1}^{m-r}u_i+\sum_{i=1}^{r}e_i+\sum_{i=1}^{r}\deg(\alpha_i)=rd.
\end{equation}

\end{theorem}

As a consequence  of  (\ref{eqIST}), the eigenstructure of a polynomial matrix determines its degree.

Given a polynomial matrix $P(s)\in \FF[s]^{m\times n}$ of $\deg(P(s))= d$,  let $e_1\leq  \dots \leq  e_r$ be the  partial multiplicities of $\infty$ in $P(s)$  and  $p_1\leq \dots \leq p_r$ the invariant orders at $\infty$  of $P(s)$.
Then (see \cite[Proposition 6.14]{AmMaZa15}), 
\begin{equation*}\label{eqqf}
e_i=p_i+d,\quad 1\leq i \leq r.
  \end{equation*}
As a consequence, $\deg (P(s))=-p_1$.
Hence, knowing the  degree and  the partial multiplicities of $\infty$ in $P(s)$ is the same as knowing its invariant orders at $\infty$; i.e., the information provided by the complete structural data is equivalent to that provided by the  eigenstructure.
Thus, we can restate Theorem \ref{theoexistenceDeDoVa152} as follows.

\begin{theorem}\label{theoexistencepolioi}
Let $m$, $n$, $r\leq \min\{m,n\}$ be  positive integers.
Let
$\alpha_1(s)\mid \dots \mid \alpha_r(s)$ be monic  polynomials.
Let $p_1\leq \dots \leq p_r$ be integers and $(c_1, \ldots, c_{n-r})$, $(u_1, \ldots, u_{m-r})$
be partitions. 
Then, there exists a polynomial matrix  $P(s)\in \FF[s]^{m\times n}$ of  $\rank(P(s))= r$,
with  $\alpha_1(s), \dots, \alpha_r(s)$ as invariant factors, $p_1, \dots, p_r$ as invariant orders at $\infty$, and $c_1, \dots, c_{n-r}$ and $u_1, \dots, u_{m-r}$ as column and row minimal indices,
 respectively, if and only if
 \begin{equation*}\label{eqISTbis}
\sum_{i=1}^{n-r}c_i+\sum_{i=1}^{m-r}u_i+\sum_{i=1}^{r}p_i+\sum_{i=1}^{r}\deg(\alpha_i)=0.
\end{equation*}
\end{theorem}

Now we state our first problem, which is a generalization of the row completion problem of polynomial matrices solved in \cite[Theorem 4.2]{AAmBaMaRo24}. 

\begin{problem}\label{problempolioi}
Let $P(s)\in\efe[s]^{m\times n}$ be a polynomial matrix. Find necessary and sufficient conditions for the existence of a polynomial  matrix  $W(s)\in \efe[s]^{z\times n}$
such that  $\begin{bmatrix}P(s)\\W(s)\end{bmatrix}$  has prescribed complete structural data.
\end{problem}

We would like to point out that the row completion problem studied in  \cite{AAmBaMaRo24} requires that $\deg \begin{bmatrix}P(s)\\W(s)\end{bmatrix} =\deg P(s)$. Here this restriction is removed, i.e.,  $\deg \begin{bmatrix}P(s)\\W(s)\end{bmatrix} \geq \deg P(s)$.

It is our aim to also study the row completion problem for rational matrices. First of all we extend Theorem \ref{theoexistencepolioi}.  
The following lemma is essential to generalize to rational matrices  some results obtained for polynomial matrices.

 \begin{lemma}\label{lem_polrat}
Let $R(s)$ be a rational matrix and let $\psi(s)$ be a monic polynomial multiple	of the least common denominator of the entries in $R(s)$. Then, $\psi(s)R(s)$ is a polynomial matrix of the same rank as $R(s)$ and
\begin{itemize}
	\item[(i)] the quotients  $\frac{\eta_1(s)}{\varphi_1(s)},\dots,\frac{\eta_r(s)}{\varphi_r(s)}$ are the invariant rational functions of $R(s)$ if and only if the polynomials $\frac{\psi(s)\eta_1(s)}{\varphi_1(s)},\dots,\frac{\psi(s)\eta_r(s)}{\varphi_r(s)}$ are the invariant factors of $\psi(s)R(s)$. 
	\item[(ii)]
 the integers $\tilde p_1,\dots,\tilde p_r$ are the invariant orders at $\infty$ of $R(s)$ if and only if the integers 
 $\tilde p_1-\deg(\psi(s)), \dots,\tilde p_r-\deg(\psi(s))$ are the invariant orders at $\infty$ of  $\psi(s)R(s)$. 
	\item[(iii)] $\mathcal{N}_r(R(s))=\mathcal{N}_r(\psi(s)R(s))$, $\mathcal{N}_{\ell}(R(s))=\mathcal{N}_{\ell}(\psi(s)R(s))$ and, therefore, the minimal indices of $\psi(s)R(s)$ and of $R(s)$ are the same.
\end{itemize}
\end{lemma}

{\bf Proof.}  Items $(i)$ and $(ii)$ can be easily derived from the  Smith--McMillan forms. The proof of item $(iii)$  is straightforward.
\hfill $\Box$

\medskip
When $\FF$ is an infinite field, Theorem 4.1 of \cite{AnDoHoMa19} provides necessary and sufficient conditions for the existence of a rational matrix $R(s)\in \FF(s)^{m\times n}$ with prescribed complete structural data. The proof is based on Theorem 3.3 of \cite{DeDoVa15}, which establishes an analogous result for polynomial matrices over infinite fields. This theorem was generalized to arbitrary fields in \cite[Theorem 3.1]{AAmBaMaRo24}.
Using the latter result,    we obtain a generalization of   \cite[Theorem 4.1]{AnDoHoMa19} to arbitrary fields.

\begin{theorem}\label{theoexistencerat}
Let $m, n, r\leq \min\{m,n\}$  be  positive integers.
Let
$\eta_1(s)\mid\dots \mid  \eta_r(s)$ and 
 $\varphi_r(s)\mid\dots\mid \varphi_1(s)$ be monic polynomials
such that  $\frac{\eta_1(s)}{\varphi_1(s)},\dots,\frac{\eta_r(s)}{\varphi_r(s)}$ are irreducible rational functions.
Let $\tilde p_1\leq \dots \leq \tilde p_r$ be  integers and $(c_1, \ldots, c_{n-r})$, $(u_1, \ldots, u_{m-r})$
 partitions. 
 Then, there exists a rational matrix $R(s)\in \FF(s)^{m\times n}$, $\rank(R(s))= r$,  with  $\frac{\eta_1(s)}{\varphi_1(s)},\dots,\frac{\eta_r(s)}{\varphi_r(s)}$ as invariant rational functions,   $\tilde p_1, \dots, \tilde p_r$  as invariant orders at $\infty$, and $c_1, \dots, c_{n-r}$ and $u_1, \dots, u_{m-r}$ as column and row minimal indices, respectively, if and only if
\begin{equation*}\label{eqISTrat}
\sum_{i=1}^{n-r}c_i+\sum_{i=1}^{m-r}u_i+\sum_{i=1}^{r}\tilde p_i+\sum_{i=1}^{r}\deg(\eta_i)-\sum_{i=1}^{r}\deg(\varphi_i)=0.
\end{equation*}
\end{theorem}
    {\bf Proof.} The proof is analogous to that of \cite[Theorem 4.1]{AnDoHoMa19} using Theorem 3.1 of \cite{AAmBaMaRo24} instead of Theorem 3.3 of \cite{DeDoVa15}, and    Lemma \ref{lem_polrat}.
\hfill $\Box$

\medskip
Now, we generalize to rational matrices Problem \ref{problempolioi}.

\begin{problem}\label{problemrat}
Let $R(s)\in\efe(s)^{m\times n}$ be a rational matrix. Find necessary and sufficient conditions for the existence of a rational matrix  $\widetilde W(s)\in \efe(s)^{z\times n}$
such that  $\begin{bmatrix}R(s)\\ \widetilde W(s)\end{bmatrix}$  has prescribed complete structural data.
\end{problem}

We are also interested in solving row (column) completion problems when we prescribe part of the complete structural data, i.e., when one or some of the four types of invariants which form the complete structural data are prescribed.

\begin{problem}\label{problempart}
Let $R(s)\in\efe(s)^{m\times n}$ $(P(s)\in\efe[s]^{m\times n})$. Find necessary and sufficient conditions for the existence of a matrix  $\widetilde W(s)\in \efe(s)^{z\times n}$ $(W(s)\in \efe[s]^{z\times n})$
such that  $\begin{bmatrix}R(s)\\ \widetilde W(s)\end{bmatrix}$ $\left(\begin{bmatrix}P(s)\\W(s)\end{bmatrix}\right)$ has part of the structural data prescribed.
\end{problem}

A solution  to Problem \ref{problempart} for polynomial matrices when only the finite structure is prescribed follows from a well-known result: the  characterization  of the invariant factors of a polynomial matrix with a  prescribed submatrix (see the next Theorem \ref{lemmaSaTh}).

\begin{theorem}[{\rm\cite[Chapter 7]{Bo64}, \cite{Sa79}, \cite{Th79}}] 
  \label{lemmaSaTh}
  Let $P(s)\in\efe[s]^{m\times n}$ and $Q(s)\in\efe[s]^{(m+z)\times (n+q)}$  be  polynomial matrices of  $\rank (P(s))=r$ and $\rank(Q(s))=\bar r$, 
  and let $\alpha_1(s)\mid \dots \mid \alpha_r(s)$ and  $\beta_1(s)\mid \dots \mid \beta_{\bar r}(s)$ be the invariant factors of $P(s)$ and $Q(s)$, respectively. There exist matrices
$X(s)\in \efe[s]^{m\times q}$, $Y(s)\in \efe[s]^{z\times q}$, $W(s)\in \efe[s]^{z\times n}$ such that $\begin{bmatrix}P(s)&X(s)\\W(s)&Y(s)\end{bmatrix}$ is unimodularly equivalent to $Q(s)$ if and only if
  $$
\beta_i(s)\mid \alpha_i(s)\mid \beta_{i+z+q}(s), \quad 1\leq i\leq r.
$$
\end{theorem}

For the  rational case and when only the finite structure is prescribed, the following result gives a solution to Problem \ref{problempart}.

\begin{theorem}[{\rm\cite[Theorem 1]{Ba89}}] 
  \label{lemmainvrat}
  Let $R(s)\in\efe(s)^{m\times n}, G(s)\in\efe(s)^{(m+z)\times (n+q)}$  be  rational matrices, $\rank (R(s))=r$, $\rank(G(s))=\bar r$.
  Let $\frac{\eta_1(s)}{\varphi_1(s)}, \dots,  \frac{\eta_r(s)}{\varphi_r(s)}$ and  $\frac{\epsilon_1(s)}{\psi_1(s)}, \dots, \frac{\epsilon_{\bar r}(s)}{\psi_{\bar r}(s)}$ be the invariant rational functions  of $R(s)$ and $G(s)$, respectively. There exist matrices
$X(s)\in \efe(s)^{m\times q}$, $Y(s)\in \efe(s)^{z\times q}$, $W(s)\in \efe(s)^{z\times n}$ such that $\begin{bmatrix}R(s)&X(s)\\W(s)&Y(s)\end{bmatrix}$ is unimodularly equivalent to $G(s)$ if and only if
  $$
\epsilon_i(s)\mid \eta_i(s)\mid \epsilon_{i+z+q}(s), 
\quad
\psi_{i+z+q}(s)\mid \varphi_i(s)\mid \psi_{i}(s), \quad 1\leq i\leq r.
$$
\end{theorem}

\section{Row (column) completion  with prescribed complete structural data}
\label{secpolcompl} 

The aim of this section is to present a solution to Problems \ref{problempolioi} and  \ref{problemrat}.

Given a polynomial matrix $P(s)$, the \textit{grade} of $P(s)$ is 
 an integer which is at
least as large as $\deg(P(s))$ (see \cite{DeDoMa14}).
We denote it by $\gr(P(s))$.

\begin{definition}\label{defFrpbform}
  Let $P(s)=P_gs^g+P_{g-1}s^{g-1}+\cdots+P_1s+P_0\in\efe[s]^{m\times n}$ be  a polynomial matrix of grade $g\geq 1$. The {\rm first Frobenius companion form of $P(s)$ with respect to  $g$} is the $(m+(g-1)n)\times gn$ pencil $C_{g,P}(s)=sX_1+Y_1$ with
$$
X_1=\begin{bmatrix}P_g&&&\\&I_n&&\\&&\ddots&\\&&&I_n\end{bmatrix}\,\text{ and }\,
Y_1=\begin{bmatrix}P_{g-1}&P_{g-2}&\cdots&P_0\\-I_n&0&\cdots&0\\&\ddots&\ddots&\vdots\\0&&-I_n&0\end{bmatrix}.
$$
\end{definition}	

When $g=\deg(P(s))$ we omit ``with respect to  $g$'' and  $C_{g,P}(s)$ is denoted by $C_P(s)$.
Notice  that when $g=1$, $C_{g,P}(s)=P(s)$.
The following lemma is a consequence of \cite[Theorems 5.3 and 4.1]{DeDoMa14}.

\begin{lemma}\label{lem_lin}
	Let $P(s)\in\efe[s]^{m\times n}$  be a polynomial matrix of  $\gr(P(s))=g\geq 1$, and let $C_{g,P}(s)$ be its first Frobenius companion form with respect to $g$.  Then,
	\begin{enumerate}
	\item\label{itemalpha}
If $\alpha_{1}(s), \dots, \alpha_{r}(s)$  are the invariant factors of $P(s)$, then the
invariant factors of $C_{g,P}(s)$ are
 $1, \stackrel{(g-1)n}{\dots}, 1, \alpha_{1}(s), \dots, \alpha_{r}(s)$. 

\item \label{itemp}
If $p_1, \dots, p_r$  are the invariant orders at $\infty$ of  $P(s)$, then
 $-1,\stackrel{(g-1)n}{\dots},$ $-1, g-1+p_1, \dots, g-1+p_r$ 
 are the invariant orders at $\infty$ of  $C_{g,P}(s)$.

	\item\label{itemc}
	If $c_1\geq\cdots\geq c_{n-r}$ are the column minimal indices of $P(s)$, then $c_1+g-1\geq\cdots\geq c_{n-r}+g-1$ are the column minimal indices of $C_{g,P}(s)$.
\item \label{itemu}
	If $u_1\geq\cdots\geq u_{m-r}$ are the row minimal indices of $P(s)$, then $u_1\geq\cdots\geq u_{m-r}$ are also the row minimal indices of $C_{g,P}(s)$.
        \end{enumerate}       
\end{lemma}

As a consequence of Lemma \ref{lem_lin}   we obtain  the next corollary.

\begin{corollary}\label{corcompofse}
Let $P(s), \bar P(s)\in\efe[s]^{m\times n}$  such that  $\gr (P(s))= \gr (\bar P(s))=g\geq 1$, and let  $C_{g,P}(s)$, $C_{g,\bar P}(s)$ be their respective first Frobenius companion forms with respect to $g$. Then,
$P(s)\approx \bar P(s)$ if and only if $C_{g,P}(s)\se C_{g, \bar P}(s)$.	
  \end{corollary}

In  Theorem 4.3 of \cite{DoSt19} (see Theorem \ref{theopencilcompletion} below)  a solution of the row completion problem for matrix pencils is given. 
We state the result for non constant pencils.  It involves the definition of the generalized majorization  (see Definition \ref{gm} below).

Let   $\bc= (c_1,  \ldots, c_x)$ and $\ba= (a_1, \ldots, a_x)$ be  two sequences of integers. It is said that   $\bc$ is \textit{majorized} by $\ba$ (denoted by $\bc \prec \ba$) if $\sum_{i=1}^k c_i \leq \sum_{i=1}^k a_i $ for $1 \leq k \leq x-1$ and $\sum_{i=1}^x  c_i =\sum_{i=1}^x a_i$ (this is an extension to sequences of integers of the definition of majorization given for partitions in \cite{HLP88}).

\begin{definition}{\rm \cite[Definition 2]{DoStEJC10}}\label{gm}
Let $\bd = (d_1, \dots, d_{q-x})$, $\ba=(a_1, \dots, a_x)$ and $\bc=(c_1, \dots, c_{q})$  be sequences of  integers.
We say that  $\bc$ is majorized by $\bd$ and $\ba$  $(\bc \prec' (\bd,\ba))$ if
\begin{equation}\label{gmaj1}
d_i\geq c_{i+x}, \quad 1\leq i\leq q-x,
\end{equation}
\begin{equation}\label{gmaj2}
\sum_{i=1}^{h_j}c_i-\sum_{i=1}^{h_j-j}d_i\leq \sum_{i=1}^j a_i, \quad 1\leq j\leq x,
\end{equation}
where $h_j=\min\{i\; : \; d_{i-j+1}<c_i\}$, $1\leq j\leq x$ 
 $(d_{q-x+1}=-\infty)$,
\begin{equation}\label{gmaj3}
\sum_{i=1}^{q}c_i=\sum_{i=1}^{q-x}d_i+\sum_{i=1}^xa_i.
\end{equation}
\end{definition}

    In the case that $x=0$, condition (\ref{gmaj2}) disappears, and conditions (\ref{gmaj1}) and (\ref{gmaj3}) are equivalent to  $\bc=\bd$. 
On the other hand, if $q=x$ then $\bc \prec' (\bd,\ba)$ is equivalent to  $\bc \prec \ba$.

\begin{theorem}[\mbox{\rm \cite[Theorem 4.3]{DoSt19}}]
{\rm (Prescription of the complete structural data for non constant matrix pencils)}
  \label{theopencilcompletion}
Let $C(s)\in  \efe[s]^{(\bar r+p)\times (\bar r+q)}$ be a matrix
pencil, $\deg(C(s))=1$, $\rank(C(s))=\bar r$. Let $\bar \phi_{1}(s,t)\mid \dots  \mid\bar \phi_{\bar r}(s,t)$  be its homogeneous invariant factors, $\bar \bc=(\bar c_1,  \dots, \bar c_q)$  its
column minimal indices, and $\bar \bu=(\bar u_1, \dots, \bar u_p)$  its row minimal
indices, where   $\bar u_1 \geq \dots \geq \bar u_{\theta}  > \bar u_{\theta +1}= \dots = \bar u_p=0 $. Let $x$ and $y$ be non negative integers. 
Let $D(s) \in  \efe[s]^{(\bar r+p+x+y)\times (\bar r+q)}$ be a matrix pencil,  $\rank(D(s))=\bar r+x$. Let $\bar \gamma _{1}(s,t) \mid \dots  \mid\bar \gamma _{\bar r+x}(s,t)$
be its homogeneous invariant factors, $\bar \bd=(\bar d_1,  \dots,  \bar d_{q - x})$  its column minimal
indices, and  $\bar  \bv=(\bar v_1, \dots, \bar v_{p+y})$   its row minimal indices, where $\bar v_1 \geq  \dots  \geq  \bar v_{\bar \theta}  > \bar v_{\bar \theta +1} = \dots =\bar v_{p+y} = 0$.
There exists a pencil $A(s)$ such that
$\begin{bmatrix}C(s)\\A(s)\end{bmatrix}\se D(s)$
if and only if
\begin{equation}\label{eqinterhom}
 \bar \gamma _{i}(s,t)\mid \bar \phi _{i}(s,t)\mid \bar \gamma _{ i+x+y}(s,t),\quad 1\leq i \leq \bar r,
\end{equation}
\begin{equation}\label{eqtheta}
\bar \theta \geq  \theta,
\end{equation}
\begin{equation}\label{eqcmimaj}
 \bar \bc  \prec'  (\bar \bd , \bar \ba),
\end{equation}
\begin{equation}\label{eqrmimaj}
\bar  \bv \prec'  (\bar \bu , \bar \bb),
\end{equation}
\begin{equation}\label{eqdegsum}
\sum_{i=1}^{\bar r+x}\deg(\lcm(\bar \phi _{i-x},\bar \gamma _{i}))\leq \sum_{i=1}^{p+y}\bar v_i-\sum_{i=1}^{p}\bar u_i+\sum_{i=1}^{\bar r+x}\deg(\bar \gamma _{ i}),
\end{equation}
where $\bar \ba = (\bar a_1, \dots, \bar a_x )$ and $\bar \bb = (\bar b_1, \dots,  \bar b_y )$ are defined as
\begin{equation*}\label{eqdefbara}
  \begin{array}{rl}
    \sum_{i=1}^{j}\bar a_i=&
\sum_{i=1}^{p+y}\bar v_i-\sum_{i=1}^{p}\bar u_i+\sum_{i=1}^{\bar r+x}\deg(\bar \gamma _{ i})\\&-\sum_{i=1}^{\bar r+x-j}\deg(\lcm(\bar\phi _{i-x+j},\bar \gamma _{ i}))-j,\quad  1\leq j \leq x,
\end{array}
\end{equation*}
\begin{equation*}\label{eqdefbarb}
  \begin{array}{rl}
    \sum_{i=1}^{j}\bar b_i
=&
\sum_{i=1}^{p+y}\bar v_i-\sum_{i=1}^{p}\bar u_i+\sum_{i=1}^{\bar r+x}\deg(\bar \gamma _{ i})\\&-\sum_{i=1}^{\bar r+x}\deg(\lcm(\bar\phi _{i-x-j},\bar \gamma _{ i})),\quad 1\leq j \leq y.
  \end{array}
\end{equation*}
\end{theorem}

\begin{rem}\label{rempenciolsioi} \

\begin{enumerate}
\item 	\label{remdecreasing1} 
Let $\bar r,x$ and $y$ be non negative integers. Given two polynomial chains
$\bar \phi_1(s,t)\mid\cdots\mid\bar \phi_{\bar r}(s,t)$ and
$\bar \gamma_1(s,t)\mid\cdots\mid\bar \gamma_{\bar r+x}(s,t)$, by \cite[Lemmas 1 and 2]{DoStEJC10}  
we can see that 
  for $1\leq j \leq x-1$,
  $$
\begin{array}{rl}
     &\sum_{i=1}^{  \bar r+x-j+1}\deg(\lcm(  \bar \phi_{i-x+j-1},
\bar \gamma_{i}))-\sum_{i=1}^{  \bar r+x-j}\deg(\lcm(  \bar \phi_{i-x+j}, \bar \gamma_{i}))
\\ \geq &\sum_{i=1}^{  \bar r+x-j}\deg(\lcm(  \bar \phi_{i-x+j},
\bar \gamma_{i}))-\sum_{i=1}^{  \bar r+x-j-1}\deg(\lcm(  \bar \phi_{i-x+j+1},
\bar \gamma_{i})),
\end{array}
$$
and for $1\leq j \leq y-1$,
$$
\begin{array}{rl}
&\sum_{i=1}^{  \bar r+x}\deg(\lcm( \bar \phi_{i-x-j+1}, \bar \gamma_i))-\sum_{i=1}^{ 
\bar r+x}\deg(\lcm( \bar \phi_{i-x-j}, \bar \gamma_i))\\\geq&
\sum_{i=1}^{  \bar r+x}\deg(\lcm( \bar \phi_{i-x-j}, \bar \gamma_i))-\sum_{i=1}^{ 
\bar r+x}\deg(\lcm( \bar \phi_{i-x-j-1}, \bar \gamma_i)).
\end{array}
$$
As a consequence, from (\ref{eqdegsum}), in Theorem \ref{theopencilcompletion}  
we obtain that $\bar a_1\geq \dots\geq   \bar a_x$ and $\bar b_1\geq \dots\geq  
\bar b_{y}\geq 0$.

Along the paper, finite sequences of integers  similar to $\bar a_1, \dots,   \bar a_x$ or $\bar b_1, \dots,  \bar b_{y},$ will be introduced. They will analogously be   decreasing. We will omit the explanation.

\item In Theorem \ref{theopencilcompletion}, let $\bar \alpha_1(s), \dots, \bar \alpha_{\bar r}(s)$ and $\bar p_1, \dots, \bar p_{\bar r}$ be the invariant factors and
the invariant orders at $\infty$ of $C(s)$, respectively, and let 
$\bar \beta_1(s), \dots, \bar \beta_{\bar r+x}(s)$ and  $\bar q_1, \dots, \bar q_{\bar r+x}$ be the invariant factors and the invariant orders at $\infty$ of  $D(s)$, respectively.
Then,
$$
\begin{array}{l}
\bar\phi_{i}(s,t)=t^{\bar p_i+1}t^{\deg(\bar\alpha_i)}\bar \alpha_i(\frac{s}{t}), \quad 1\leq i \leq \bar r,\\
\bar\gamma_{i}(s,t)=t^{\bar q_i+1}t^{\deg(\bar\beta_i)}\bar \beta_i(\frac{s}{t}), \quad 1\leq i \leq \bar r+x.
\end{array}
$$
Hence, (\ref{eqinterhom}) is equivalent to 
\begin{equation}\label{eqinterifpencils}
       \bar \beta_i(s)\mid \bar \alpha_{i}(s)\mid \bar\beta_{i+x+y}(s),\quad 1\leq i \leq \bar r,
       \end{equation}
       \begin{equation}\label{eqinterioipencils}
       \quad \bar q_i \leq \bar p_i \leq \bar q_{i+x+y}, 
       \quad 1\leq i \leq \bar r,
 \end{equation}
 and (\ref{eqdegsum}) is equivalent to
\begin{equation}\label{eqdegsumpolpencilsioi}
   \begin{array}{rl}
  &\sum_{i=1}^{\bar r}\deg(\lcm(\bar \alpha_i,\bar \beta_{i+x}))+
  \sum_{i=1}^{\bar r}\max\{\bar p_i,\bar q_{i+x}\}\\
  \leq& \sum_{i=1}^{p+y}\bar v_i-\sum_{i=1}^{p}\bar u_i
  +\sum_{i=1}^{\bar r}\deg(\bar \beta_{i+x}) +\sum_{i=1}^{\bar r}\bar q_{i+x}.\end{array}
 \end{equation}
	Moreover,
 \begin{equation*}\label{eqdefapencilsioi}
		\begin{array}{rl}
			  \sum_{i=1}^{j}\bar a_i=&
			\sum_{i=1}^{p+y}\bar v_i-\sum_{i=1}^{p} \bar u_i+\
   \sum_{i=1}^{\bar r+j}\deg( \bar \beta_{i+x-j})+\sum_{i=1}^{\bar r+j}\bar q_{i+x-j}
                        \\&
   -\sum_{i=1}^{\bar r}\deg(\lcm(\bar \alpha_{i},  \bar \beta_{i+x-j})
                        -\sum_{i=1}^{\bar r}\max\{\bar p_i,\bar q_{i+x-j}\},\\&\hfill
			1\leq j \leq x,
		\end{array}
\end{equation*}
\begin{equation*}\label{eqdefbpencilsioi}
		\begin{array}{rl}
\sum_{i=1}^{j}\bar  b_i=&
			\sum_{i=1}^{p+y}  \bar v_i-\sum_{i=1}^{p}  \bar u_i
                        +\sum_{i=1}^{  \bar r-j}\deg(  \bar \beta_{i+x+j})+\sum_{i=1}^{  \bar r-j}
                        \bar q_{i+x+j}\\&
			-\sum_{i=1}^{  \bar r-j}\deg(\lcm(\bar \alpha_i,  \bar \beta_{i+x+j}))-
\sum_{i=1}^{  \bar r-j}\max\{\bar p_i,  \bar q_{i+x+j}\},
 \\&\hfill 1 \leq j \leq y.
		\end{array}
	\end{equation*} 

\end{enumerate}
\end{rem}

\begin{proposition}\label{propeqprob}
  Let $P(s)\in\efe[s]^{m\times n}$ and $Q(s)\in\efe[s]^{(m+z)\times n}$ be such that   $\deg(Q(s))=g\ge \max\{\deg(P(s)),1\}$. Let $C_{g,P}(s)$ be the first Frobenius companion form of $P(s)$ with respect to $g$ and $C_{Q}(s)$ be the first Frobenius companion form of $Q(s)$.
  Then, there exists $W(s)\in\efe[s]^{z\times n}$ such that  $\begin{bmatrix}P(s)\\W(s)\end{bmatrix}\approx Q(s)$ if and only if there exists a matrix pencil
   $A(s)\in\efe[s]^{z\times gn}$ such that $\begin{bmatrix}C_{g,P}(s)\\A(s)\end{bmatrix}\se C_Q(s)$.
  \end{proposition}
{\bf Proof.}
The proof is completely analogous to that of  \cite[Proposition 4.1]{AAmBaMaRo24} 
exchanging degree and grade and applying  Corollary \ref{corcompofse}.
\hfill $\Box$

\medskip
Now, we can give a solution to Problem \ref{problempolioi}.

\begin{theorem}{\rm (Prescription of the complete structural data for  polynomial matrices)}
  \label{theoprescr4ioi}
  Let $P(s)\in\efe[s]^{m\times n}$ be a polynomial matrix of  $\rank (P(s))=r$.
	Let $\alpha_1(s)\mid \cdots\mid\alpha_{r}(s)$ be its 
invariant factors, $p_1, \dots, p_r$ its invariant orders at $\infty$,
	$\bc=(c_1,  \dots,  c_{n-r})$  its
	column minimal indices, and $\bu=(u_1, \dots,  u_{m-r})$  its row minimal,  where   $u_1 \geq \dots \geq u_{\eta}  >  u_{\eta +1}= \dots = u_{m-r}=0 $.
	 
  Let $z,x$ be integers such that $0\leq x\leq \min\{z, n-r\}$. 
	 	Let $\beta_1(s)\mid \cdots\mid\beta_{r+x}(s)$ be monic polynomials, $q_1\leq\dots\leq q_{r+x}$ integers, and 
	$\bd=(d_1,  \dots,  d_{n-r-x})$  and $\bv=(v_1, \dots,  v_{m+z-r-x})$ two partitions, where
	$v_1 \geq  \dots  \geq  v_{\bar \eta}  > v_{\bar \eta +1} = \dots =v_{m+z-r-x} = 0$.
	There exists a polynomial matrix  $W(s)\in \efe[s]^{z\times n}$ such that  	
 $\rank \left(\begin{bmatrix}P(s)\\W(s)\end{bmatrix}\right) =r+x$ and 
	$\begin{bmatrix}P(s)\\W(s)\end{bmatrix}$ has $\beta_1(s), \dots, \beta_{r+x}(s)$ as invariant factors,
	$q_1, \dots, q_{r+x}$ as invariant orders at $\infty$,
 $d_1, \dots, d_{n-r-x}$ as column minimal indices
	and 
	$v_1, \dots, v_{m+z-r-x}$ as row minimal indices  
	if and only if
\begin{equation}\label{eqinterif}
        \beta_i(s)\mid \alpha_{i}(s)\mid \beta_{i+z}(s),\quad 1\leq i \leq r,
 \end{equation}
\begin{equation}\label{eqinterpolioi}
q_i\leq p_i\leq q_{i+z}, \quad 1\leq i\leq r,
        \end{equation}
        \begin{equation}\label{eqetapol}\bar \eta \geq  \eta,\end{equation}
	\begin{equation}\label{eqcmimajpol}  \bc \prec'  (\bd , \ba),\end{equation}
	\begin{equation}\label{eqrmimajpol}\bv \prec'  (\bu , \bb),\end{equation}
	\begin{equation}\label{eqdegsumpolioi}
   \begin{aligned}
  \sum_{i=1}^{r}\deg(\lcm(\alpha_i,\beta_{i+x}))+
  \sum_{i=1}^{r}\max\{p_i,q_{i+x}\}\\
  \leq \sum_{i=1}^{m+z-r-x}v_i-\sum_{i=1}^{m-r}u_i
  +\sum_{i=1}^{r}\deg(\beta_{i+x}) +\sum_{i=1}^{r}q_{i+x},\\
 \mbox{ with equality when $x=0$,}
\end{aligned}
 \end{equation}
	where $ \ba = (a_1, \dots, a_x )$ and $\bb = ( b_1, \dots,  b_{z-x} )$ are defined as
	\begin{equation}\label{eqdefaioi}
 \begin{array}{rl}
			  \sum_{i=1}^{j}a_i=&
			\sum_{i=1}^{m+z-r-x}v_i-\sum_{i=1}^{m-r} u_i+
   \sum_{i=1}^{ r+j}\deg(  \beta_{i+x-j})+\sum_{i=1}^{ r+j}q_{i+x-j}
                        \\&
   -\sum_{i=1}^{ r}\deg(\lcm( \alpha_{i},  \beta_{i+x-j}))
                        -\sum_{i=1}^{r}\max\{ p_i, q_{i+x-j}\},\\&\hfill
			1\leq j \leq x,
		\end{array}
\end{equation}
\begin{equation}\label{eqdefbioi}
		\begin{array}{rl}
\sum_{i=1}^{j} b_i=&
\sum_{i=1}^{m+z-r-x}v_i-\sum_{i=1}^{m-r} u_i
                        +\sum_{i=1}^{r-j}\deg(\beta_{i+x+j})+\sum_{i=1}^{r-j}
                        q_{i+x+j}\\&
			-\sum_{i=1}^{r-j}\deg(\lcm(\alpha_i, \beta_{i+x+j}))-
\sum_{i=1}^{r-j}\max\{ p_i,  q_{i+x+j}\},\\& \hfill
  1 \leq j \leq z-x.
		\end{array}
  \end{equation}

\end{theorem}

{\bf Proof.} 
The proof is analogous to that of \cite[Theorem 4.2]{AAmBaMaRo24}.
Define $d=-p_1$ and $g=-q_1$.
Then $\deg(P(s))=d$. 

When  $g\geq d$,  we can build $C_{g,P}(s)$, the first Frobenius companion form of $P(s)$ with respect to $g$. If 
 $g\geq 1$, we will take 
$$
\bar r=(g-1)n+r, \; y=z-x, \; p=m-r=m+(g-1)n-\bar r, \; q=n-r=gn-\bar r.
$$

Assume that there exists  a polynomial matrix  $W(s)\in \efe[s]^{z\times n}$ such that  	
	$Q(s)=\begin{bmatrix}P(s)\\W(s)\end{bmatrix}$ has the prescribed invariants. Then, $\deg(Q(s))=g\geq d$.
If $g=0$, then
        (\ref{eqinterif})-(\ref{eqdegsumpolioi}) trivially hold.
If $g\geq 1$, 
let $C_{g,Q}(s)=C_Q(s)$ be the first Frobenius companion form of $Q(s)$. 
 By Proposition \ref{propeqprob}, there exists a matrix pencil  $A(s)\in \efe[s]^{z\times gn}$  such that $\begin{bmatrix}C_{g,P}(s)\\A(s)\end{bmatrix}\se C_{Q}(s)$.
 
 Let $\bar \alpha_{1}(s),\cdots,\bar \alpha_{\bar r}(s)$,
 $\bar p_1,\dots, \bar p_{\bar r}$,
 $\bar \bc= (\bar c_1,  \dots, \bar c_{q})$ and
 $\bar \bu= (\bar u_1,  \dots, \bar u_{p})$
 be the invariant factors,
invariant orders at $\infty$,
column minimal indices and row minimal indices of $C_{g,P}(s)$,
where $\bar u_1 \geq  \dots  \geq  \bar u_{\theta}  > \bar u_{\theta +1} = \dots =\bar u_{p} = 0$ and let
$\bar \beta_{1}(s),\cdots,\bar \beta_{\bar r+x}(s)$,
 $\bar q_1,\dots, \bar q_{\bar r+x}$,
$\bar \bd= (\bar d_1,  \dots, \bar d_{q-x})$ and
 $\bar \bv= (\bar v_1,  \dots, \bar v_{p+y})$
be the invariant factors,
invariant orders at $\infty$,
column minimal indices and row minimal indices of $C_{Q}(s)$,
where
 $\bar v_1 \geq  \dots  \geq  \bar v_{\bar \theta}  > \bar v_{\bar \theta +1} = \dots =\bar v_{p+y} = 0$.
By  Theorem \ref{theopencilcompletion} and Remark \ref{rempenciolsioi}, conditions (\ref{eqtheta})--(\ref{eqrmimaj}) and 
(\ref{eqinterifpencils})--(\ref{eqdegsumpolpencilsioi}) hold.
Applying  Lemma \ref{lem_lin}, it is easy to see that 
  (\ref{eqtheta})--(\ref{eqrmimaj}) and 
(\ref{eqinterifpencils})--(\ref{eqdegsumpolpencilsioi}) are equivalent to (\ref{eqinterif})--(\ref{eqdegsumpolioi}).

Assume now that 
(\ref{eqinterif})--(\ref{eqdegsumpolioi}) hold.
Then, from (\ref{eqcmimajpol}), (\ref{eqdegsumpolioi}) if $x=0$, (\ref{eqinterif}) and (\ref{eqinterpolioi}), we get 
$\sum_{i=1}^{n-r} c_i-\sum_{i=1}^{n-r-x} d_i=\sum_{i=1}^xa_i=\sum_{i=1}^{m+z-r-x}v_i-\sum_{i=1}^{m-r} u_i
                        +\sum_{i=1}^{r+x}\deg( \beta_{i})+
                        \sum_{i=1}^{r+x}q_{i}-
                         \sum_{i=1}^{r}\deg( \alpha_{i})
                        -\sum_{i=1}^{r}p_{i}
$.  By Theorem \ref{theoexistencepolioi} applied to $P(s)$, we obtain
$$0=\sum_{i=1}^{n-r-x} d_i
+\sum_{i=1}^{m+z-r-x}v_i
                        +\sum_{i=1}^{r+x}\deg( \beta_{i})+
                        \sum_{i=1}^{r+x}q_{i}                        
.$$
Applying again Theorem \ref{theoexistencepolioi}, we derive that there exists 
a polynomial matrix  $Q(s)\in \FF[s]^{(m+z)\times n}$, $\rank(Q(s))= r+x$,
with  $\beta_1(s), \dots, \beta_{r+x}(s)$ as invariant factors, $q_1, \dots, q_{r+x}$ as invariant orders at $\infty$, and $d_1, \dots, d_{n-r-x}$ and $v_1, \dots, v_{m+z-r-x}$ as column and row minimal indices. Then, $\deg(Q(s))=-q_1=g$.  From (\ref{eqinterpolioi}) we obtain $g=-q_1\geq -p_1=d$.

If $g=0$, then choosing $W\in \efe^{z\times n}$  such that 
	$\rank \begin{bmatrix}P\\W\end{bmatrix}=r+x$, the matrix $\begin{bmatrix}P\\W\end{bmatrix}$ has the prescribed invariants.
    If $g\geq 1$,
let $C_Q(s)$  be the first Frobenius companion form of $Q(s)$ and let $\bar \beta_1(s), \dots, \bar \beta_{\bar r+x}(s)$, $\bar q_1, \dots, \bar q_{\bar r+x}$,
$\bar \bd=(\bar d_1, \dots, \bar d_{q-x})$ and  $\bar \bv=(\bar v_1,  \dots, \bar v_{p+y})$ be the invariant factors, invariant orders at $\infty$, column minimal indices and row minimal indices of $C_Q(s)$, respectively,
where $\bar v_1 \geq  \dots  \geq  \bar v_{\bar \eta} > \bar v_{\bar \eta+1} = \dots =\bar v_{p+y} = 0$.
As in the proof of the necessity, (\ref{eqinterif})--(\ref{eqdegsumpolioi}) are equivalent to 
(\ref{eqtheta})--(\ref{eqrmimaj}) and 
(\ref{eqinterifpencils})--(\ref{eqdegsumpolpencilsioi}).
The result follows from
Theorem \ref{theopencilcompletion}, Remark \ref{rempenciolsioi} and 
 Proposition \ref{propeqprob}.  
\hfill $\Box$

\begin{rem}\label{remcomparison}
Under the  conditions of Theorem \ref{theoprescr4ioi}, let 
$e_i=p_i-p_1$, $1\leq i \leq r$, and 
$f_i=q_i-q_1$, $1\leq i \leq r+x$. 
Then, conditions (\ref{eqinterpolioi}) and  (\ref{eqdegsumpolioi}) become
\begin{equation*}\label{eqinterpolif}
f_i\leq e_i+p_1-q_1\leq f_{i+z}, \quad 1\leq i\leq r,       \end{equation*}
  and      
	\begin{equation*}\label{eqdegsumpolif}
   \begin{aligned}
  \sum_{i=1}^{r}\deg(\lcm(\alpha_i,\beta_{i+x}))+
  \sum_{i=1}^{r}\max\{e_i+p_1-q_1,f_{i+x}\}\\
  \leq \sum_{i=1}^{m+z-r-x}v_i-\sum_{i=1}^{m-r}u_i
  +\sum_{i=1}^{r}\deg(\beta_{i+x}) +\sum_{i=1}^{r}f_{i+x},\\
 \mbox{ with equality when $x=0$,}
\end{aligned}
 \end{equation*}
 respectively, and
	$ \ba = (a_1, \dots, a_x )$ and $\bb = ( b_1, \dots,  b_{z-x} )$ can be rewritten as
	\begin{equation*}\label{eqdefaif}
 \begin{array}{rl}
			  \sum_{i=1}^{j}a_i=&
			\sum_{i=1}^{m+z-r-x}v_i-\sum_{i=1}^{m-r} u_i+
   \sum_{i=1}^{ r+j}\deg(  \beta_{i+x-j})+\sum_{i=1}^{ r+j}f_{i+x-j}
                        \\&
   -\sum_{i=1}^{ r}\deg(\lcm( \alpha_{i},  \beta_{i+x-j}))
                        -\sum_{i=1}^{r}\max\{ e_i+p_1-q_1, f_{i+x-j}\}\\&+jq_1, \quad
			1\leq j \leq x,
		\end{array}
\end{equation*}
\begin{equation*}\label{eqdefbif}
		\begin{array}{rl}
\sum_{i=1}^{j} b_i=&
\sum_{i=1}^{m+z-r-x}v_i-\sum_{i=1}^{m-r} u_i
                        +\sum_{i=1}^{r-j}\deg(\beta_{i+x+j})+\sum_{i=1}^{r-j}
                        f_{i+x+j}\\&
			-\sum_{i=1}^{r-j}\deg(\lcm(\alpha_i, \beta_{i+x+j}))-
\sum_{i=1}^{r-j}\max\{ e_i+p_1-q_1,  f_{i+x+j}\},\\& \hfill
  1 \leq j \leq z-x.
		\end{array}
  \end{equation*}
Therefore, when $q_1=p_1=-\deg(P(s))$, from Theorem \ref{theoprescr4ioi} we recover  \cite[Theorem 4.2]{AAmBaMaRo24}.
\end{rem}

 In order to solve Problem \ref{problemrat} we will use the following technical lemma.

\begin{lemma}\label{lemma_lcm}
Let $\varphi(s),\eta(s),\psi(s),\epsilon(s),\pi(s)\in\efe[s]$ such that $\varphi(s)\mid \pi(s)$, $\psi(s)\mid \pi(s)$, $\gcd(\varphi,\eta)=1$ and $\gcd(\psi,\epsilon)=1$. Then,
$$
\lcm\left(\frac{\pi}{\varphi}\eta, \frac{\pi}{\psi}\epsilon\right)=
\frac{\pi(s)}{\gcd(\varphi,\psi)}\lcm(\eta,\epsilon).
$$
\end{lemma}
{\bf Proof.} Note that  $h_1(s)=\frac{\varphi(s)}{\gcd(\varphi,\psi)}\frac{\lcm(\eta,\epsilon)}{\eta(s)}$ and 
$h_2(s)=\frac{\psi(s)}{\gcd(\varphi,\psi)}\frac{\lcm(\eta,\epsilon)}{\epsilon(s)}$ are polynomials. Then,
$$\frac{\pi(s)}{\varphi(s)}\eta(s) h_1(s)=\frac{\pi(s)}{\gcd(\varphi,\psi)}\lcm(\eta,\epsilon),\quad
\frac{\pi(s)}{\psi(s)}\epsilon(s) h_2(s)=\frac{\pi(s)}{\gcd(\varphi,\psi)}\lcm(\eta,\epsilon),
$$
and both $\frac{\pi(s)}{\varphi(s)}\eta(s)$ and $\frac{\pi(s)}{\psi(s)}\epsilon(s)$ are divisors of $\frac{\pi(s)}{\gcd(\varphi,\psi)}\lcm(\eta,\epsilon)$. 
Therefore, $\lcm\left(\frac{\pi}{\varphi}\eta, \frac{\pi}{\psi}\epsilon\right)$ divides $\frac{\pi(s)}{\gcd(\varphi,\psi)}\lcm(\eta,\epsilon)$, i.e., there exists $q(s)\in\efe[s]$ such that 
$$\lcm\left(\frac{\pi}{\varphi}\eta, \frac{\pi}{\psi}\epsilon\right)=\frac{\pi(s)}{\gcd(\varphi,\psi)}\lcm(\eta,\epsilon)\frac{1}{q(s)}.$$ Let $\ell_1(s), \ell_2(s)\in\efe[s]$ such that
$$
\begin{array}{c}
\frac{\pi(s)}{\varphi(s)}\eta(s)\ell_1(s)=\frac{\pi(s)}{\gcd(\varphi,\psi)}\lcm(\eta,\epsilon)\frac{1}{q(s)},\\
\frac{\pi(s)}{\psi(s)}\epsilon(s)\ell_2(s)=\frac{\pi(s)}{\gcd(\varphi,\psi)}\lcm(\eta,\epsilon)\frac{1}{q(s)}.
\end{array}
$$
Then, $q(s)$ is a divisor of $h_1(s)$ and of $h_2(s)$. 

Let $x_1(s)=\gcd(\eta, \epsilon)$. Then $\eta(s)=x_1(s)x_2(s)$ and $\epsilon(s)=x_1(s)x_3(s)$ with  $\gcd(x_2,x_3)=1$. Thus, $\lcm(\eta,\epsilon)=x_1(s)x_2(s)x_3(s)$, $\frac{\lcm(\eta,\epsilon)}{\eta(s)}=x_3(s)$, and $\frac{\lcm(\eta,\epsilon)}{\epsilon(s)}=x_2(s)$.
Note that 
$h_1(s)=\frac{\varphi(s)}{\gcd(\varphi,\psi)}x_3(s)$ and  $h_2(s)=\frac{\psi(s)}{\gcd(\varphi,\psi)}x_2(s)$. 

Let $q(s)=q_1(s)q_2(s)$ with $q_1(s), q_2(s)\in\efe[s]$ such that $q_1(s)$ divides $\frac{\varphi(s)}{\gcd(\varphi,\psi)}$ and $q_2(s)$ divides $x_3(s)$. As $\gcd(\varphi,\eta)=1$,  $\gcd(q_1,\eta)=1$ and $\gcd(q_1,x_2)=1$. Thus, $q_1(s)$ divides $\frac{\psi(s)}{\gcd(\varphi,\psi)}$. 
Analogously, since $\gcd(x_2,x_3)=1$, $\gcd(q_2,x_2)=1$ and  $q_2(s)$ divides $\frac{\psi(s)}{\gcd(\varphi,\psi)}$. It follows from $\gcd(\psi,\epsilon)=1$ that $\gcd(q_2,\epsilon)=1$ and $\gcd(q_2,x_3)=1$. Thus, $q_2(s)=1$, and $q(s)=q_1(s)$
is a divisor of both $\frac{\varphi(s)}{\gcd(\varphi,\psi)}$ and $\frac{\psi(s)}{\gcd(\varphi,\psi)}$. Hence, $q(s)=1$ and the result follows. 
\hfill $\Box$

\medskip

In the sequel we use the following notation:
given $\varphi(s),\eta(s),\psi(s),\epsilon(s)\in\efe[s]$ such that  $\gcd(\varphi,\eta)=1$ and $\gcd(\psi,\epsilon)=1$, and $p, q$ integers, we denote
$$
\Delta\left(\frac{\eta}{\varphi}, \frac{\epsilon}{\psi}, p, q\right)=
\deg(\lcm(\eta,\epsilon))- \deg(\gcd(\varphi,\psi)) +\max\{p,q\},
$$
$$\Delta\left(\frac{\eta}{\varphi}, \frac{\epsilon}{\psi}\right)=
\deg(\lcm(\eta,\epsilon))- \deg(\gcd(\varphi,\psi)),$$
$$\Delta\left(\frac{\eta}{\varphi}, p\right)=\deg(\eta)- \deg(\varphi)+p,$$

$$\Delta\left(\frac{\eta}{\varphi}\right)=
\deg(\eta)- \deg(\varphi).$$

\begin{theorem}{\rm (Prescription of the complete structural data for rational matrices)}\label{prescr4rat}
  Let $R(s)\in\efe(s)^{m\times n}$ be a rational matrix, $\rank (R(s))=r$. 
 Let $\frac{\eta_1(s)}{\varphi_1(s)},\dots, \frac{\eta_r(s)}{\varphi_r(s)}$ 
	be its invariant rational functions,  
$\tilde p_1, \dots, \tilde p_r$ its invariant orders at $\infty$, 
	$\bc=(c_1,  \dots,  c_{n-r})$  its
	column minimal indices, and $\bu=(u_1, \dots,  u_{m-r})$  its row minimal
	indices, where   $u_1 \geq \dots \geq u_{\eta}  >  u_{\eta +1}= \dots = u_{m-r}=0 $.
	
	Let $z,x$ be integers such that $0\leq x\leq \min\{z, n-r\}$ and let $\epsilon_1(s)\mid \dots \mid\epsilon_{r+x}(s)$ and 
        $\psi_{r+x}(s)\mid \dots \mid\psi_{1}(s)$ be monic polynomials such that $\frac{\epsilon_i(s)}{\psi_i(s)}$ are irreducible rational functions, $1\leq i \leq r+x$.
        Let $\tilde q_{1}\leq \dots \leq \tilde q_{r+x}$ be integers and $\bd=(d_1,  \dots,  d_{n-r-x})$  and $\bv=(v_1, \dots,  v_{m+z-r-x})$  be two partitions, where   
	$v_1 \geq  \dots  \geq  v_{\bar \eta}  > v_{\bar \eta +1} = \dots =v_{m+z-r-x} = 0$.
	There exists a rational matrix  $\widetilde W(s)\in \efe(s)^{z\times n}$ such that $\rank \left(\begin{bmatrix}R(s)\\\widetilde W(s)\end{bmatrix}\right)=r+x$ and
        $\begin{bmatrix}R(s)\\\widetilde W(s)\end{bmatrix}$ has $\frac{\epsilon_1(s)}{\psi_1(s)}, \dots, \frac{\epsilon_{r+x}(s)}{\psi_{r+x}(s)}$ as invariant rational functions, 
 $\tilde q_1, \dots,\tilde q_{r+x}$ as invariant orders at $\infty$,  $d_1,  \dots,  d_{n-r-x}$  as column minimal indices and $v_1, \dots,  v_{m+z-r-x}$ as row minimal indices  
        if and only if 
(\ref{eqetapol}),
        \begin{equation}\label{eqinterratnum}
		\epsilon_i(s)\mid \eta_i(s)\mid \epsilon_{i+z}(s),\quad 1\leq i \leq r,
	\end{equation}
        \begin{equation}\label{eqinterratden}
		\psi_{i+z}(s)\mid \varphi_i(s)\mid \psi_{i}(s),\quad 1\leq i \leq r,
	\end{equation}
        \begin{equation}\label{eqinterratioi}
\tilde q_i\leq \tilde p_i\leq \tilde q_{i+z}, \quad 1\leq i\leq r,
        \end{equation}
	\begin{equation}\label{eqcmimajrat}  \bc \prec'  (\bd , \tilde \ba),\end{equation}
	\begin{equation}\label{eqrmimajrat}\bv \prec'  (\bu , \tilde \bb),\end{equation}
\begin{equation}\label{eqdegsumrat}
\begin{aligned}
  \sum_{i=1}^{r}
  \Delta\left(\frac{\eta_i}{\varphi_i}, \frac{\epsilon_{i+x}}{\psi_{i+x}}, \tilde p_i,\tilde q_{i+x}\right)-\sum_{i=1}^{r}\Delta\left(\frac{\epsilon_{i+x}}{\psi_{i+x}}, \tilde q_{i+x}\right)
  \leq \sum_{i=1}^{m+z-r-x}v_i-\sum_{i=1}^{m-r}u_i
  ,\\
  \mbox{ with equality when $x=0$,}
\end{aligned}
 \end{equation}
	where $\tilde \ba = (\tilde a_1, \dots, \tilde a_x )$ and $\tilde \bb = (\tilde b_1, \dots,  \tilde b_{z-x} )$ are defined as
 \begin{equation}\label{eqdeftildea}
		\begin{array}{rl}
		\sum_{i=1}^{j}	\tilde a_i=&
			\sum_{i=1}^{m+z-r-x}v_i-\sum_{i=1}^{m-r} u_i+\sum_{i=1}^{r+j}\Delta\left(\frac{\epsilon_{i+x-j}}{\psi_{i+x-j}}, \tilde q_{i+x-j}\right)\\&
-\sum_{i=1}^{  r} \Delta\left(\frac{\eta_i}{\varphi_i}, \frac{\epsilon_{i+x-j}}{\psi_{i+x-j}}, \tilde p_{i}, \tilde q_{i+x-j}\right),\quad
                        1\leq j \leq x,
		\end{array}
\end{equation}
\begin{equation}\label{eqdeftildeb}
	\begin{array}{rl}
 \sum_{i=1}^{j}\tilde b_i=&
			\sum_{i=1}^{m+z-r-x}  v_i-\sum_{i=1}^{m-r}  u_i+\sum_{i=1}^{  r-j}\Delta\left(\frac{\epsilon_{i+x+j}}{\psi_{i+x+j}},\tilde q_{i+x+j}\right)\\&
-\sum_{i=1}^{  r-j}
\Delta\left(\frac{\eta_i}{\varphi_i}, \frac{\epsilon_{i+x+j}}{\psi_{i+x+j}}, \tilde p_i,  \tilde q_{i+x+j}\right),\quad
                        1\leq j \leq z-x,\\
		\end{array}
	\end{equation}  
\end{theorem}
{\bf Proof.}  We start with a remark  assuming that  $\varphi_1(s) \mid \psi_1(s)$.
Define
$
d=\deg (\psi_1)-\tilde p_1$,  $g=\deg (\psi_1)-\tilde q_1$,
$$
\begin{array}{l}
\alpha_i(s)=\psi_1(s)\frac{\eta_i(s)}{\varphi_i(s)},
\quad p_i=\tilde p_i-\deg(\psi_1),
\quad 1\leq i \leq r,\\
\beta_i(s)=\psi_1(s)\frac{\epsilon_i(s)}{\psi_i(s)},\quad 
q_i=\tilde q_i-\deg(\psi_1), \quad  1\leq i \leq r+x.\\
\end{array}
$$
Then
(\ref{eqinterif}) is equivalent to
$$
\epsilon_i(s)\varphi_i(s)\mid \eta_i(s)\psi_i(s), \quad
\eta_i(s)\psi_{i+z}(s)\mid \epsilon_{i+z}(s)\varphi_i(s), \quad 1\leq i \leq r.
$$
As $\gcd(\epsilon_i, \psi_i)=1$, $1\leq i \leq r+x$ and
$\gcd(\eta_i, \varphi_i)=1$, $1\leq i \leq r$, we derive that (\ref{eqinterif}) is equivalent to (\ref{eqinterratnum}) and (\ref{eqinterratden}).
It is clear that (\ref{eqinterpolioi}) is equivalent to (\ref{eqinterratioi}).
Define also $\ba=(a_1, \dots, a_x)$ and $\bb=(b_1, \dots, b_{z-x})$ as in
(\ref{eqdefaioi}) and (\ref{eqdefbioi}).
Then, by Lemma \ref{lemma_lcm}, 
$\ba=\tilde \ba$ and $\bb=\tilde \bb$, hence (\ref{eqcmimajpol}) and 
	(\ref{eqrmimajpol})  are equivalent to (\ref{eqcmimajrat}) and 
	(\ref{eqrmimajrat}),  respectively. Analogously,  
condition (\ref{eqdegsumpolioi}) is equivalent to
(\ref{eqdegsumrat}).

Assume that there exists a
rational matrix  $\widetilde W(s)\in \efe(s)^{z\times n}$ such that, if
$G(s)=\begin{bmatrix}R(s)\\\widetilde W(s)\end{bmatrix}$, then $G(s)$ has the prescribed structural data.
Recall that $\psi_1(s)$ and $\varphi_1(s)$  are the monic least common denominator of the entries of $G(s)$ and $R(s)$, respectively.
Thus, the matrix $\psi_1(s)G(s)=\begin{bmatrix}\psi_1(s)R(s)\\\psi_1(s)\widetilde{W}(s)\end{bmatrix}$ is polynomial and $\varphi_1(s)\mid \psi_1(s)$.
Let $$P(s)=\psi_1(s)R(s), \quad Q(s)=\psi_1(s)G(s).$$
By Lemma \ref{lem_polrat} we  know that  $\rank(P(s))=r$, 
$\alpha_1(s), \dots, \alpha_r(s)$ are the invariant factors,
$p_1, \dots, p_r$ the invariant orders at $\infty$,
$c_1, \dots, c_{n-r}$ the  column minimal indices
and $u_1, \dots, u_{m-r}$ the  row minimal indices of $P(s)$, and $\rank(Q(s))=r+x$,
 $\beta_1(s), \dots, \beta_{r+x}(s)$ are the invariant factors,
$q_1, \dots, q_{r+x}$ the invariant orders at $\infty$,
 $d_1, \dots, d_{n-r-x}$ the  column 
 and $v_1, \dots, v_{m+z-r-x}$ the  row minimal indices
 of $Q(s)$.
By Theorem \ref{theoprescr4ioi},
(\ref{eqinterif})--(\ref{eqdegsumpolioi})
hold, where $\ba$ and $\bb$ are defined in 
(\ref{eqdefaioi}) and (\ref{eqdefbioi}), respectively. Equivalently, (\ref{eqetapol}) and 
(\ref{eqinterratnum})-(\ref{eqdegsumrat})
hold, where $\tilde \ba$ and $\tilde \bb$ are defined as in (\ref{eqdeftildea}) and  (\ref{eqdeftildeb}), respectively.

Conversely, assume that (\ref{eqetapol}) and 
(\ref{eqinterratnum})--(\ref{eqdegsumrat}) are satisfied. Then 
$\varphi_1(s)\mid \psi_1(s)$ and 
(\ref{eqinterif})--(\ref{eqdegsumpolioi}) hold.
Let $P(s)=\psi_1(s)R(s)$.
By Lemma \ref{lem_polrat}, $\rank(P(s))=r$,
$\alpha_1(s), \dots, \alpha_r(s)$ are the invariant factors,
$p_1, \dots, p_r$ the invariant orders at $\infty$,
$c_1, \dots, c_{n-r}$ the  column minimal indices
and $u_1, \dots, u_{m-r}$ the  row minimal indices of $P(s)$.

From (\ref{eqinterif})--(\ref{eqdegsumpolioi})  by Theorem \ref{theoprescr4ioi} 
there exists a polynomial matrix  $W(s)\in \efe[s]^{z\times n}$ such that  	
$\rank \left(\begin{bmatrix}P(s)\\W(s)\end{bmatrix}\right)=r+x$ and
        $\begin{bmatrix}P(s)\\W(s)\end{bmatrix}$ has $\beta_1(s), \dots, \beta_{r+x}(s)$ as invariant factors, 
 $q_1, \dots, q_{r+x}$ as invariant orders at $\infty$,  $d_1,  \dots,  d_{n-r-x}$  as column minimal indices and $v_1, \dots,  v_{m+z-r-x}$ as row minimal indices.
Let $\widetilde W(s)=\frac{1}{\psi_1(s)}W(s)$. Then  $\widetilde W(s)\in \efe(s)^{z\times n}$ and $\begin{bmatrix}P(s)\\W(s)\end{bmatrix}=
\psi_1(s)\begin{bmatrix}R(s)\\\widetilde W(s)\end{bmatrix}
$.
By Lemma \ref{lem_polrat}, $\rank \left(\begin{bmatrix}R(s)\\\widetilde W(s)\end{bmatrix}\right)=r+x$ and
        $\begin{bmatrix}R(s)\\\widetilde W(s)\end{bmatrix}$ has $\frac{\epsilon_1(s)}{\psi_1(s)}, \dots, \frac{\epsilon_{r+x}(s)}{\psi_{r+x}(s)}$ as invariant rational functions, 
 $\tilde q_1, \dots, \tilde q_{r+x}$ as invariant orders at $\infty$,  $d_1,  \dots,  d_{n-r-x}$  as column minimal indices and $v_1, \dots,  v_{m+z-r-x}$ as row minimal indices.
\hfill $\Box$

\begin{rem}\label{remratcdioi}
By  Theorem \ref{theoexistencerat}  (see also \cite[Remark 4.3]{AAmBaMaRo24}),
if (\ref{eqetapol}) and (\ref{eqinterratnum})--(\ref{eqdegsumrat})
 hold, then
 \begin{equation}\label{eqdegsumratcdioi}
\begin{aligned}
  \sum_{i=1}^{r}\Delta\left(\frac{\eta_i}{\varphi_i}, \frac{\epsilon_{i+x}}{\psi_{i+x}}, \tilde p_i,\tilde q_{i+x}\right) 
  \leq \sum_{i=1}^{n-r}c_i-\sum_{i=1}^{n-r-x}d_i
  +\sum_{i=1}^{r}\Delta\left(\frac{\eta_i}{\varphi_i}, \tilde p_i\right) 
  -  \sum_{i=1}^{x}\Delta\left(\frac{\epsilon_{i}}{\psi_{i}},\tilde q_{i}\right),\\
\mbox{ with equality when $x=z$,} 
  \end{aligned}
 \end{equation}
and (\ref{eqcmimajrat}) and  (\ref{eqrmimajrat}) hold for
   $\tilde \ba=(\tilde a_1, \dots, \tilde a_x)$ and $\tilde \bb=(\tilde b_1, \dots, \tilde b_{z-x})$ defined as
   \begin{equation}\label{eqdeftildeacdioi}
                \begin{array}{rl}
			 \sum_{i=1}^{j}\tilde a_i=& \sum_{i=1}^{n-r}c_i-\sum_{i=1}^{n-r-x}d_i+
                        \sum_{i=1}^{r}\Delta\left(\frac{\eta_i}{\varphi_i}, \tilde p_i\right) 
  -  \sum_{i=1}^{x-j}\Delta\left(\frac{\epsilon_{i}}{\psi_{i}},\tilde q_{i}\right)
  \\&-\sum_{i=1}^{r}\Delta\left(\frac{\eta_i}{\varphi_i}, \frac{\epsilon_{i+x-j}}{\psi_{i+x-j}}, \tilde p_i,\tilde q_{i+x-j}\right),
               \quad 1\leq j \leq x,
		\end{array}
\end{equation}
\begin{equation}\label{eqdeftildebcdbioi}
	\begin{array}{rl}
  \sum_{i=1}^{j}\tilde b_i=&
\sum_{i=1}^{n-r}c_i-\sum_{i=1}^{n-r-x}d_i
+
                        \sum_{i=1}^{r}\Delta\left(\frac{\eta_i}{\varphi_i}, \tilde p_i\right)-\sum_{i=1}^{x+j}\Delta\left(\frac{\epsilon_{i}}{\psi_{i}},\tilde q_{i}\right) \\
   &-\sum_{i=1}^{r-j}\Delta\left(\frac{\eta_i}{\varphi_i}, \frac{\epsilon_{i+x+j}}{\psi_{i+x+j}}, \tilde p_i,\tilde q_{i+x+j}\right),
                        \quad
                        1\leq j \leq z-x.
		\end{array}
	\end{equation}
Conversely,
(\ref{eqetapol}), (\ref{eqinterratnum})--(\ref{eqrmimajrat}) and 
(\ref{eqdegsumratcdioi}) 
with $\tilde \ba$ and $\tilde \bb$   defined as in 
(\ref{eqdeftildeacdioi}) and (\ref{eqdeftildebcdbioi}), respectively, imply
(\ref{eqetapol}) and (\ref{eqinterratnum})--(\ref{eqdegsumrat}) 
 with $\tilde \ba$ and $\tilde \bb$   defined as in 
(\ref{eqdeftildea}) and (\ref{eqdeftildeb}), respectively.  

  \end{rem}

\section{Row (column) completion with part of the structural data prescribed}
\label{secpart}

In this section we first solve Problem \ref{problempart} when the complete structural data  but the row (column) minimal indices are prescribed (see Subsection \ref{subsec_fininfcolrow}). Afterwards, in Subsection \ref{subsec_dinoronf}, we solve Problem \ref{problempart} 
when the finite and/or infinite structures are prescribed.

Given a rational matrix $R(s)\in\efe(s)^{m\times n}$ with  
 $\frac{\eta_1(s)}{\varphi_1(s)},\dots, \frac{\eta_r(s)}{\varphi_r(s)}$ 
as  invariant rational functions, when $\varphi_1(s)=1$ the matrix $R(s)\in\efe[s]^{m\times n}$ is polynomial with invariant factors $\eta_1(s),\dots, \eta_r(s)$.
When we prescribe the invariant rational functions (Theorems \ref{theoprescrifioicmi}, \ref{theoprescrifioirmi}, \ref{theoprescrifioi} 
and \ref{theoprescrifdegree})   
we  present the results for rational matrices, and the polynomial cases are derived from them.

\subsection{Prescription of the finite and infinite structures and column or row minimal indices}\label{subsec_fininfcolrow}

We present two results related to Problem \ref{problempart}. In Theorem  \ref{theoprescrifioicmi} we prescribe the finite and infinite structures and column  minimal indices, and in Theorem \ref{theoprescrifioirmi} we replace the column minimal indices by the row minimal indices. The proofs are analogous to those of \cite[Sections 4.2, 4.3]{AAmBaMaRo24}.

\begin{theorem}
\label{theoprescrifioicmi}
{\rm (Prescription of the finite and infinite structures, and the column minimal indices)}
  Let $R(s)\in\efe(s)^{m\times n}$ be a rational matrix, $\rank (R(s))=r$. 
 Let $\frac{\eta_1(s)}{\varphi_1(s)},\dots, \frac{\eta_r(s)}{\varphi_r(s)}$ 
	be its invariant rational functions,  
$\tilde p_1, \dots, \tilde p_r$ its invariant orders at $\infty$ and
	$\bc=(c_1,  \dots,  c_{n-r})$  its
	column minimal indices. 
	
	Let $z,x$ be integers such that $0\leq x\leq \min\{z, n-r\}$ and let $\epsilon_1(s)\mid \dots \mid\epsilon_{r+x}(s)$ and 
        $\psi_{r+x}(s)\mid \dots \mid\psi_{1}(s)$ be monic polynomials such that $\frac{\epsilon_i(s)}{\psi_i(s)}$ are irreducible rational functions, $1\leq i \leq r+x$.
        Let $\tilde q_{1}\leq \dots \leq \tilde q_{r+x}$ be integers and $\bd=(d_1,  \dots,  d_{n-r-x})$  a partition. 
	There exists a rational matrix  $\widetilde W(s)\in \efe(s)^{z\times n}$ such that $\rank \left(\begin{bmatrix}R(s)\\\widetilde W(s)\end{bmatrix}\right)=r+x$ and
        $\begin{bmatrix}R(s)\\\widetilde W(s)\end{bmatrix}$ has $\frac{\epsilon_1(s)}{\psi_1(s)}, \dots, \frac{\epsilon_{r+x}(s)}{\psi_{r+x}(s)}$ as invariant rational functions, 
 $\tilde q_1, \dots, \tilde q_{r+x}$ as invariant orders at $\infty$ and  $d_1,  \dots,  d_{n-r-x}$  as column minimal indices
        if and only if
         (\ref{eqinterratnum})--(\ref{eqcmimajrat})  and (\ref{eqdegsumratcdioi}), 
where $\tilde \ba = (\tilde a_1, \dots, \tilde a_x )$ is  defined as in (\ref{eqdeftildeacdioi}).
  \end{theorem}
{\bf Proof.} It is analogous to the  proof  of Theorem 4.5 of \cite{AAmBaMaRo24}.
\hfill $\Box$

\begin{theorem} {\rm (Prescription of the finite and infinite structures, and the row minimal indices)}
  \label{theoprescrifioirmi}
 Let $R(s)\in\efe(s)^{m\times n}$ be a rational matrix, $\rank (R(s))=r$. 
 Let $\frac{\eta_1(s)}{\varphi_1(s)},\dots, \frac{\eta_r(s)}{\varphi_r(s)}$ 
	be its invariant rational functions,  
$\tilde p_1, \dots, \tilde p_r$ its invariant orders at $\infty$, 
	$\bc=(c_1,  \dots,  c_{n-r})$  its
	column minimal indices, and $\bu=(u_1, \dots,  u_{m-r})$  its row minimal
	indices, where   $u_1 \geq \dots \geq u_{\eta}  >  u_{\eta +1}= \dots = u_{m-r}=0 $.
	
	Let $z,x$ be integers such that $0\leq x\leq \min\{z, n-r\}$ and let $\epsilon_1(s)\mid \dots \mid\epsilon_{r+x}(s)$ and 
        $\psi_{r+x}(s)\mid \dots \mid\psi_{1}(s)$ be monic polynomials such that $\frac{\epsilon_i(s)}{\psi_i(s)}$ are irreducible rational functions, $1\leq i \leq r+x$.
        Let $\tilde q_{1}\leq \dots \leq \tilde q_{r+x}$ be integers and $\bv=(v_1,  \dots,  v_{m+z-r-x})$  a partition  such that
        $v_1 \geq  \dots  \geq  v_{\bar \eta} > v_{\bar \eta+1}=\dots= v_{m+z-r-x}=0$.
Let
$\tilde \ba = (\tilde a_1, \dots, \tilde a_x )$ and $\tilde \bb = (\tilde b_1, \dots,  \tilde b_{z-x} )$
be defined as in (\ref{eqdeftildea}) and (\ref{eqdeftildeb}), respectively.

\begin{enumerate}
\item
  If $x=n-r$, 
        there exists a rational matrix  $\widetilde W(s)\in \efe(s)^{z\times n}$ such that $\rank \left(\begin{bmatrix}R(s)\\\widetilde W(s)\end{bmatrix}\right)=r+x$ and
        $\begin{bmatrix}R(s)\\\widetilde W(s)\end{bmatrix}$ has $\frac{\epsilon_1(s)}{\psi_1(s)}, \dots, \frac{\epsilon_{r+x}(s)}{\psi_{r+x}(s)}$ as invariant rational functions, 
 $\tilde q_1, \dots, \tilde q_{r+x}$ as invariant orders at $\infty$ and  $v_1,  \dots,  v_{m+z-r-x}$  as row minimal indices
        if and only if
        (\ref{eqetapol}),
        (\ref{eqinterratnum})--(\ref{eqinterratioi}), (\ref{eqrmimajrat}), (\ref{eqdegsumrat}) and
    \begin{equation*}\label{eqcprectildea}\bc \prec \tilde\ba.\end{equation*}  
\item
  If $x<n-r$, 
        there exists a rational matrix  $\widetilde W(s)\in \efe(s)^{z\times n}$ such that $\rank \left(\begin{bmatrix}R(s)\\\widetilde W(s)\end{bmatrix}\right)=r+x$ and
        $\begin{bmatrix}R(s)\\\widetilde W(s)\end{bmatrix}$ has $\frac{\epsilon_1(s)}{\psi_1(s)}, \dots, \frac{\epsilon_{r+x}(s)}{\psi_{r+x}(s)}$ as invariant rational functions, 
 $\tilde q_1, \dots, \tilde q_{r+x}$ as invariant orders at $\infty$ and  $v_1,  \dots,  v_{m+z-r-x}$  as row minimal indices
        if and only if
        (\ref{eqetapol}),
        (\ref{eqinterratnum})--(\ref{eqinterratioi}), (\ref{eqrmimajrat}), (\ref{eqdegsumrat}),
        \begin{equation*}\label{eqsumca2}
\sum_{i=1}^{x+1}c_i-c_\ell\geq \sum_{i=1}^x\tilde a_i,
        \end{equation*}
        and
\begin{equation*}\label{eqsumca3}
\sum_{i=j+2}^{x+1}c_i\geq \sum_{i=j+1}^x\tilde a_i,\quad \ell\leq j \leq x-1,
\end{equation*}
where $\ell=\min\{j\geq 1\; : \; \sum_{i=1}^{j}c_i>\sum_{i=1}^{j}\tilde a_i\}$.
\end{enumerate}
  \end{theorem}
{\bf Proof.} It is analogous to the  proof  of Theorem 4.8 of \cite{AAmBaMaRo24}.
\hfill $\Box$

\subsection{Prescription of the finite and/or infinite structures}\label{subsec_dinoronf}

In this subsection we deal with Problem \ref{problempart} when only the finite or the infinite structures are prescribed. First,  
we present  a solution  when both the finite and infinite structures are prescribed.  
Secondly, we give a solution when only the infinite structure is prescribed (in both polynomial and rational cases). 
As mentioned, the solution for the case where only the finite structure  is prescribed is known (for the polynomial case when the degree is not prescribed see Theorem \ref{lemmaSaTh}, and  for the rational case see Theorem \ref{lemmainvrat}). In Theorem \ref{theoprescrifdegree} we prescribe the finite structure and the first invariant order at infinity. Note that in the polynomial case prescribing the first invariant order at infinity is the same as prescribing the degree.
 
Although some of the proofs in this subsection are analogous to those presented in \cite[Section 4.4]{AAmBaMaRo24}, we write them for the convenience of the reader. 

\begin{theorem} {\rm (Prescription of the finite and infinite structures)}
  \label{theoprescrifioi}
 Let $R(s)\in\efe(s)^{m\times n}$ be a rational matrix, $\rank (R(s))=r$. 
 Let $\frac{\eta_1(s)}{\varphi_1(s)},\dots, \frac{\eta_r(s)}{\varphi_r(s)}$ 
	be its invariant rational functions,  
$\tilde p_1, \dots, \tilde p_r$ its invariant orders at $\infty$, 
	$\bc=(c_1,  \dots,  c_{n-r})$  its
	column minimal indices, and $\bu=(u_1, \dots,  u_{m-r})$  its row minimal
	indices.
	
	Let $z,x$ be integers such that $0\leq x\leq \min\{z, n-r\}$, let $\epsilon_1(s)\mid \dots \mid\epsilon_{r+x}(s)$ and 
        $\psi_{r+x}(s)\mid \dots \mid\psi_{1}(s)$ be monic polynomials such that $\frac{\epsilon_i(s)}{\psi_i(s)}$ are irreducible rational functions, $1\leq i \leq r+x$, and let 
        $\tilde q_{1}\leq \dots \leq \tilde q_{r+x}$ be integers. 
  \begin{enumerate}
  \item  \label{ittheoprescrifio11} If  $x<z$ or $x=z=n-r$, then
there exists a
 rational matrix  $\widetilde W(s)\in \efe(s)^{z\times n}$ such that $\rank \left(\begin{bmatrix}R(s)\\\widetilde W(s)\end{bmatrix}\right)=r+x$ and
        $\begin{bmatrix}R(s)\\\widetilde W(s)\end{bmatrix}$ has $\frac{\epsilon_1(s)}{\psi_1(s)}, \dots, \frac{\epsilon_{r+x}(s)}{\psi_{r+x}(s)}$ as invariant rational functions, 
 $\tilde q_1, \dots,  \tilde q_{r+x}$ as invariant orders at $\infty$
        if and only if
        (\ref{eqinterratnum})--(\ref{eqinterratioi}) and 
        \begin{equation}\label{eqfromhifj}
   \begin{aligned}
            \sum_{i=1}^{r}\Delta\left(\frac{\eta_i}{\varphi_i}, \frac{\epsilon_{i+x-j}}{\psi_{i+x-j}},\tilde p_i, \tilde q_{i+x-j}\right)
            +\sum_{i=1}^{x-j}\Delta\left(\frac{\epsilon_{i}}{\psi_{i}},\tilde q_{i}\right)\\
            +\sum_{i=1}^{m-r}u_i+\sum_{i=1}^{j}c_i+\sum_{i=x+1}^{n-r}c_i\leq  0,
 \quad  0\leq j \leq x-1,\\
 \mbox{with equality for $j=0$ when $x=z=n-r$.}
  \end{aligned}
  \end{equation}
\item \label{ittheoprescrifioi2}
  If  $x=z<n-r$, then  there exists  $\widetilde W(s)\in \efe(s)^{z\times n}$ such that $\rank \left(\begin{bmatrix}R(s)\\\widetilde W(s)\end{bmatrix}\right)=r+x$ and
        $\begin{bmatrix}R(s)\\\widetilde W(s)\end{bmatrix}$ has $\frac{\epsilon_1(s)}{\psi_1(s)}, \dots, \frac{\epsilon_{r+x}(s)}{\psi_{r+x}(s)}$ as invariant rational functions, 
 $\tilde q_1, \dots,  \tilde q_{r+x}$ as invariant orders at $\infty$
        if and only if
        (\ref{eqinterratnum})--(\ref{eqinterratioi}),
        \begin{equation}\label{eqsumcap2}
\sum_{i=1}^{x+1}c_i-c_\ell\geq \sum_{i=1}^x\tilde a'_i,
        \end{equation}
        and
\begin{equation}\label{eqsumcap3}
\sum_{i=j+2}^{x+1}c_i\geq \sum_{i=j+1}^x\tilde a'_i,\quad \ell\leq j \leq x-1,
\end{equation}
where 
$\tilde \ba'=(\tilde a'_1, \dots, \tilde a'_x)$ is defined as
\begin{equation}\label{eqdeftildeap}
\begin{array}{l}
		\sum_{i=1}^{j}	\tilde a'_i=
		\sum_{i=1}^{r+j}\Delta\left(\frac{\epsilon_{i+x-j}}{\psi_{i+x-j}}, \tilde q_{i+x-j}\right)
-\sum_{i=1}^{  r} \Delta\left(\frac{\eta_i}{\varphi_i}, \frac{\epsilon_{i+x-j}}{\psi_{i+x-j}}, \tilde p_{i}, \tilde q_{i+x-j}\right), \\\hfill  
                        1\leq j \leq x,\end{array}
\end{equation}
and
$\ell=\min\{j\geq 1\; : \; \sum_{i=1}^{j}c_i>\sum_{i=1}^{j}\tilde a'_i\}$.


\end{enumerate}

        \end{theorem}
{\bf Proof.} The proof follows the scheme of   that of Theorem 4.10 of \cite{AAmBaMaRo24}, but it deserves some hints.

\begin{enumerate}
  \item   Case  $x<z$ or $x=z=n-r$.
Assume that there exists 
 a rational matrix  $\widetilde W(s)\in \efe(s)^{z\times n}$ 
 such that 
        $\begin{bmatrix}R(s)\\\widetilde W(s)\end{bmatrix}$ has the prescribed invariants.
Let 
$\bd=(d_1,\dots, d_{n-r-x})$ be the 
column minimal indices of $\begin{bmatrix}R(s)\\\widetilde W(s)\end{bmatrix}$.
By Theorem \ref{theoprescrifioicmi} and Remark \ref{remratcdioi},
         (\ref{eqinterratnum})--(\ref{eqcmimajrat})  and (\ref{eqdegsumratcdioi}) hold, 
where $\tilde \ba = (\tilde a_1, \dots, \tilde a_x )$ is  defined as in (\ref{eqdeftildeacdioi}).
From (\ref{eqdegsumratcdioi}) and Theorem \ref{theoexistencerat} we obtain
$$
           \sum_{i=1}^{r} \Delta\left(\frac{\eta_i}{\varphi_i}, \frac{\epsilon_{i+x}}{\psi_{i+x}}, \tilde p_i, \tilde q_{i+x}\right)+\sum_{i=1}^{x}
           \Delta\left(\frac{\epsilon_{i}}{\psi_{i}}, \tilde q_{i}\right)
            +\sum_{i=1}^{m-r}u_i+\sum_{i=1}^{n-r-x}d_i\leq  0,
$$
with equality if $x=z$.
 From (\ref{eqcmimajrat})  we get $\sum_{i=1}^{n-r-x}d_i\geq \sum_{i=1}^{n-r-x}c_{i+x}=\sum_{i=x+1}^{n-r}c_i$. Therefore,  (\ref{eqfromhifj}) holds for $j=0$.

For $1\leq j\leq x -1$, from (\ref{eqcmimajrat}), \cite[Lemma 4.9]{AAmBaMaRo24},    (\ref{eqdeftildeacdioi}) and Theorem \ref{theoexistencerat},  we obtain
$$
\begin{array}{rl}
\sum_{i=1}^{j}c_i\leq & \sum_{i=1}^{j}\tilde a_i+\sum_{i=1}^{n-r-x}d_i-\sum_{i=x+1}^{n-r}c_{i}\\
=&\sum_{i=1}^{n-r-x}d_i-\sum_{i=x+1}^{n-r}c_{i}
-\sum_{i=1}^{n-r-x} d_i
-\sum_{i=1}^{m-r}u_i\\&-\sum_{i=1}^{  r}
\Delta\left(\frac{\eta_i}{\varphi_i}, \frac{\epsilon_{i+x-j}}{\psi_{i+x-j}}, \tilde p_i,  \tilde q_{i+x-j}\right)- \sum_{i=1}^{x-j} \Delta\left(\frac{\epsilon_{i}}{\psi_{i}}, \tilde q_{i}\right).
\end{array}
$$
Thus, (\ref{eqfromhifj}) holds.

Conversely, 
assume that (\ref{eqinterratnum})--(\ref{eqinterratioi}) and (\ref{eqfromhifj})
hold. 
Define
 $\hat a_1, \dots, \hat a_x$ as 
$$
 \begin{array}{rl}
			 \sum_{i=1}^{j}\hat a_i=& \sum_{i=1}^{x}c_i+
                        \sum_{i=1}^{r}\Delta\left(\frac{\eta_i}{\varphi_i}, \tilde p_i\right) 
  -  \sum_{i=1}^{x-j}\Delta\left(\frac{\epsilon_{i}}{\psi_{i}},\tilde q_{i}\right)
  \\&-\sum_{i=1}^{r}\Delta\left(\frac{\eta_i}{\varphi_i}, \frac{\epsilon_{i+x-j}}{\psi_{i+x-j}}, \tilde p_i,\tilde q_{i+x-j}\right),
               \quad 1\leq j \leq x.
		\end{array}
$$
By condition (\ref{eqfromhifj}) for $j=0$ and Theorem \ref{theoexistencerat},
$$
            \sum_{i=1}^{r}
\Delta\left(\frac{\eta_i}{\varphi_i}, \frac{\epsilon_{i+x}}{\psi_{i+x}},\tilde p_i, \tilde q_{i+x}\right)
+\sum_{i=1}^{x}
\Delta\left(\frac{\epsilon_{i}}{\psi_{i}},\tilde q_{i}\right)
\leq
  \sum_{i=1}^{x}c_i +\sum_{i=1}^{r}
  \Delta\left(\frac{\eta_i}{\varphi_i}, \tilde p_i\right);
$$
hence
$$
\hat a_1\geq 
\sum_{i=1}^{r}
\Delta\left(\frac{\eta_i}{\varphi_i}, \frac{\epsilon_{i+x}}{\psi_{i+x}},
\tilde p_i, \tilde q_{i+x}
\right)
 -\sum_{i=1}^{  r}
\Delta\left(\frac{\eta_i}{\varphi_i}, \frac{\epsilon_{i+x-1}}{\psi_{i+x-1}}, \tilde p_i, \tilde q_{i+x-1}\right)+
\Delta\left(\frac{\epsilon_{x}}{\psi_{x}},\tilde q_{x}
\right).
$$
Taking into account Remark \ref{rempenciolsioi}.\ref{remdecreasing1}, 
 we have  
$\hat a_1\geq \hat a_2\geq\dots \geq \hat a_{x}$. 
Let  $\hat\ba=(\hat a_1, \dots, \hat a_x)$. By Theorem \ref{theoexistencerat},
for $1\leq j \leq x$,
$$
 \begin{array}{rl}
			 \sum_{i=1}^{j}\hat a_i=& -\sum_{i=x+1}^{n-r}c_i-
    \sum_{i=1}^{m-r}u_i                      
  -  \sum_{i=1}^{x-j}\Delta\left(\frac{\epsilon_{i}}{\psi_{i}},\tilde q_{i}\right)
  \\&-\sum_{i=1}^{r}\Delta\left(\frac{\eta_i}{\varphi_i}, \frac{\epsilon_{i+x-j}}{\psi_{i+x-j}}, \tilde p_i,\tilde q_{i+x-j}\right),
               \quad 1\leq j \leq x.
		\end{array}
$$
From (\ref{eqfromhifj})  we obtain 
$$\sum_{i=1}^{j}\hat a_i\geq \sum_{i=1}^{ j}   c_i,\quad 1\leq j \leq x-1.$$
Moreover, from (\ref{eqinterratnum})--(\ref{eqinterratioi}),  we obtain
$\sum_{i=1}^{ x}\hat a_i=\sum_{i=1}^{ x}c_i$.

If $n=r+x$, then $\bc \prec \hat \ba$, and let $\bd=\emptyset$ so that $\bc\prec'(\bd, \hat  \ba)$ holds. Otherwise, if $n>r+x$, by  \cite[Lemma 4.6]{AAmBaMaRo24} 
 there exists
a sequence of integers
$\bd=(d_1, \ldots, d_{n-r-x})$ such that $\bc\prec'(\bd, \hat  \ba)$, 
$d_i= c_{i+x}$ for $2\leq i \leq n-r-x$, and $d_1=\sum_{i=1}^{x+1}c_i-\sum_{i=1}^{x}\hat a_i=c_{x+1}$.

Let  $\tilde \ba=(\tilde a_1, \dots, \tilde a_{x})$
be defined as in (\ref{eqdeftildeacdioi}). 
 Then $\tilde \ba=\hat \ba$,  therefore  (\ref{eqcmimajrat}) holds.
 Furthermore, from (\ref{eqfromhifj}) for $j=0$  and Theorem \ref{theoexistencerat} we obtain
 $$
\begin{array}{ll}
  & \sum_{i=1}^{r}
\Delta\left(\frac{\eta_i}{\varphi_i}, \frac{\epsilon_{i+x}}{\psi_{i+x}},\tilde p_{i}, \tilde q_{i+x}\right)\\
\leq &
\sum_{i=1}^{n-r}(c_i-c_{i+x})+\sum_{i=1}^{r}
\Delta\left(\frac{\eta_i}{\varphi_i}, \tilde p_{i} \right)
  -\sum_{i=1}^{x}
   \Delta\left( \frac{\epsilon_{i}}{\psi_{i}}, \tilde q_{i}\right) \\= & \sum_{i=1}^{n-r}c_i-\sum_{i=1}^{n-r-x}d_i
 +\sum_{i=1}^{r}
\Delta\left(\frac{\eta_i}{\varphi_i}, \tilde p_{i} \right)
 -\sum_{i=1}^{x}
   \Delta\left( \frac{\epsilon_{i}}{\psi_{i}}, \tilde q_{i}\right),
 \end{array}
$$
 with equality if $x=z=n-r$, i.e., (\ref{eqdegsumratcdioi}) is satisfied.
 By Theorem \ref{theoprescrifioicmi},  the result follows.

\item 
  Case $x=z<n-r$.   As $x=z$,  observe that if there exits $\widetilde{W}(s)\in \efe(s)^{z\times n}$ such that $\rank\left(\begin{bmatrix}R(s)\\\widetilde{W}(s)\end{bmatrix}\right)=r+x$,  then  the row minimal indices of $\begin{bmatrix}R(s)\\\widetilde{W}(s)\end{bmatrix}$  are the row minimal indices of $R(s)$, i.e.,  $\bv=\bu$. For the sufficiency we prescribe $\bv=\bu$.  The result follows from Theorem  \ref{theoprescrifioirmi}.
\end{enumerate}
\hfill $\Box$

\begin{rem}\label{4.1819-4.17}
If $x=z<n-r$, conditions
(\ref{eqinterratnum})--(\ref{eqinterratioi}), (\ref{eqsumcap2}) and  (\ref{eqsumcap3})
imply  (\ref{eqfromhifj})
 (see \cite[Remark 4.11]{AAmBaMaRo24}). 
\end{rem}

When we only prescribe the infinite structure we present two results, one for polynomial matrices and another one for rational matrices.

\begin{theorem} {\rm(Prescription of 
 the infinite structure for polynomial matrices)}
  \label{theoprescrioipol}
 Let $P(s)\in\efe[s]^{m\times n}$ be a polynomial matrix, $\rank (P(s))=r$. 
 Let   
$p_1, \dots,  p_r$ be its invariant orders at $\infty$ and 
	$\bc=(c_1,  \dots,  c_{n-r})$  its
	column minimal indices. 
	
	Let $z,x$ be integers such that $0\leq x\leq \min\{z, n-r\}$ and let $ q_{1}\leq \dots \leq q_{r+x}$ be integers. There exists a polynomial  matrix  $W(s)\in \efe[s]^{z\times n}$ such that $\rank \left(\begin{bmatrix}P(s)\\ W(s)\end{bmatrix}\right)=r+x$ and
        $\begin{bmatrix}P(s)\\ W(s)\end{bmatrix}$ has 
 $q_1, \dots,   q_{r+x}$ as invariant orders at $\infty$
        if and only if (\ref{eqinterpolioi}) and 
\begin{equation}\label{eqfromhifjioipol}
           \sum_{i=1}^{r}\max\{ p_i,  q_{i+x-j}\}+\sum_{i=1}^{x-j} q_i-\sum_{i=1}^{r} p_i
          \leq\sum_{i=j+1}^{x}c_i,
 \quad  0\leq j \leq x-1.
        \end{equation}
        \end{theorem}
{\bf Proof.} 
The proof is similar to that of Theorem 4.12 of \cite{AAmBaMaRo24}; we precise some calculations.

Let   $\alpha_1(s),\dots, \alpha_r(s)$ 
	be the invariant  factors, and $\bu=(u_1, \dots,  u_{m-r})$  the row minimal indices of $P(s)$.

 Assume that there is a  polynomial  matrix  $W(s)\in \efe[s]^{z\times n}$, $\rank \left(\begin{bmatrix}P(s)\\ W(s)\end{bmatrix}\right)=r+x$,  such that
        $\begin{bmatrix}P(s)\\ W(s)\end{bmatrix}$ has 
 $q_1, \dots,   q_{r+x}$ as invariant orders at $\infty$.
Let  $\beta_1(s),\dots, \beta_{r+x}(s)$ 
	be its invariant factors. By Theorem \ref{theoprescrifioi} and Remark \ref{4.1819-4.17} we obtain (\ref{eqinterif}),  (\ref{eqinterpolioi}) and 
        \begin{equation}\label{eqfromhifjioipol1}
          \begin{array}{l}
  \sum_{i=1}^{r}\deg(\lcm(\alpha_{i},\beta_{i+x-j}))+ \sum_{i=1}^{x-j}\deg(\beta_i)+\sum_{i=1}^{r}\max\{ p_i,  q_{i+x-j}\}\\+\sum_{i=1}^{x-j} q_i
  +\sum_{i=1}^{m-r}u_i+\sum_{i=1}^{j}c_i+\sum_{i=x+1}^{n-r}c_i\leq  0,
  \quad  0\leq j \leq x-1.
  \end{array}
        \end{equation}
We have $$\sum_{i=1}^{r}\deg(\lcm(\alpha_{i},\beta_{i+x-j}))\geq \sum_{i=1}^{r}\deg(\alpha_i)\geq \sum_{i=1}^{r}\deg(\alpha_i)-\sum_{i=1}^{x-j}\deg(\beta_i).$$ Thus, from  (\ref{eqfromhifjioipol1}) we obtain
$$\begin{array}{l}
            \sum_{i=1}^{r}\deg(\alpha_i)+\sum_{i=1}^{r}\max\{ p_i,  q_{i+x-j}\}+\sum_{i=1}^{x-j} q_i\\
            +\sum_{i=1}^{m-r}u_i+\sum_{i=1}^{j}c_i+\sum_{i=x+1}^{n-r}c_i\leq  0,
 \quad  0\leq j \leq x-1,
  \end{array}
  $$
  which by Theorem \ref{theoexistencepolioi} is equivalent to (\ref{eqfromhifjioipol}).

Conversely, assume that 
(\ref{eqinterpolioi})
and  
(\ref{eqfromhifjioipol}) hold.
Let
$$t=\sum_{i=1}^{r}p_i+\sum_{i=1}^{x}c_i-\sum_{i=1}^{r}\max\{p_i,q_{i+x}\}-\sum_{i=1}^{x}q_i.$$
If $x=0$, from 
(\ref{eqinterpolioi})
we obtain $t=0$, and if $x>0$, from (\ref{eqfromhifjioipol}) we have $t\geq 0$.

 Define
 $$
\begin{array}{ll}
\beta_i(s)=1, &
  1\leq i \leq x,\\
  \beta_{i+x}(s)=\alpha_{i}(s), &  1\leq i \leq r-1,\\
  \beta_{r+x}(s)=\alpha_{r}(s)\tau(s),
\end{array}$$ where 
$\tau(s)$ is a monic polynomial of 
$\deg(\tau)=t$.
We have $\beta_1(s) \mid  \dots  \mid   \beta_{r+x}(s)$, and (\ref{eqinterif}) holds.

If $x>1$, for $1\leq j \leq x$, we have $\beta_{i+x-j}(s) \mid \beta_{i+x-1}(s)\mid \alpha_i(s)$, $1\leq i \leq r$, therefore
$$\sum_{i=1}^{r}\deg(\lcm(\alpha_i,\beta_{i+x-j}))+\sum_{i=1}^{x-j}\deg(\beta_i)=\left\{\begin{array}{c}
\sum_{i=1}^{r}\deg(\alpha_i)+t, \quad  j=0,\\
\sum_{i=1}^{r}\deg(\alpha_i), \quad 1\leq j  \leq x.
\end{array}\right. 
$$

Thus, from Theorem  \ref{theoexistencepolioi} and  (\ref{eqfromhifjioipol}) we obtain
\begin{equation}\label{newcondition}
\begin{array}{l}
\sum_{i=1}^{r}\deg(\lcm(\alpha_{i},\beta_{i+x-j}))+ \sum_{i=1}^{x-j}\deg(\beta_i)+\sum_{i=1}^{r}\max\{ p_i,  q_{i+x-j}\}+\sum_{i=1}^{x-j} q_i\\
  +\sum_{i=1}^{m-r}u_i+\sum_{i=1}^{j}c_i+\sum_{i=x+1}^{n-r}c_i
\\
=\left\{\begin{array}{c}
t-\sum_{i=1}^{r} p_i-\sum_{i=1}^{x}c_i   +\sum_{i=1}^{r}\max\{ p_i,  q_{i+x}\}+\sum_{i=1}^{x} q_i=0, \ j=0,\\
\sum_{i=1}^{r}\max\{ p_i,  q_{i+x-j}\}+\sum_{i=1}^{x-j} q_i-\sum_{i=1}^{r} p_i-\sum_{i=j+1}^{x} c_i \leq 0, \ 1\leq j  \leq x-1.
\end{array}\right.
\end{array}
\end{equation}

If $x<z$ or $x=z=n-r$ the result follows from  Theorem \ref{theoprescrifioi} (item \ref{ittheoprescrifio11}).

If $x=z<n-r$, 
let $\tilde \ba'=(\tilde a'_1, \dots, \tilde a'_x)$ be defined as
\begin{equation*}\label{eqdefap}
		\begin{array}{rl}
		\sum_{i=1}^{j}	\tilde a'_i=&
		\sum_{i=1}^{r+j}\deg(\beta_{i+x-j})-\sum_{i=1}^{  r} \deg(\lcm(\alpha_i, \beta_{i+x-j})\\&+\sum_{i=1}^{r+j} q_{i+x-j}
-\sum_{i=1}^{  r} \max\{p_{i},  q_{i+x-j}\},\quad
                        1\leq j \leq x.
		\end{array}
\end{equation*}
From (\ref{eqinterif}), (\ref{eqinterpolioi}) and 
(\ref{newcondition}) we have
$$\begin{array}{rl}
\sum_{i=1}^{x}	\tilde a'_i=&\sum_{i=1}^{r+x}\deg(\beta_{i})-\sum_{i=1}^{  r} \deg(\alpha_i)+\sum_{i=1}^{r+x} q_{i}
-\sum_{i=1}^{  r} p_i\\=&
\sum_{i=1}^{x}c_i-\sum_{i=1}^{r}\max\{p_i,q_{i+x}\}+\sum_{i=1}^{r}q_{i+x}= \sum_{i=1}^{x}c_i.
\end{array}
$$
Let $j\in\{1,\dots, x-1\}$. Then 
$$
\begin{array}{rl}
	\sum_{i=1}^{j}	\tilde a'_i=&
  \sum_{i=1}^{r+x}\deg(\beta_i)-\sum_{i=1}^{r}\deg(\lcm(\alpha_{i},\beta_{i+x-j}))-\sum_{i=1}^{x-j}\deg(\beta_i)\\&
  +\sum_{i=1}^{r+x}q_i-\sum_{i=1}^{r}\max\{p_{i},q_{i+x-j}\}
  -\sum_{i=1}^{x-j}q_i\\=&\sum_{i=1}^{r}\deg(\alpha_i)+\sum_{i=1}^{r}p_i+\sum_{i=1}^{x}c_i
   -\sum_{i=1}^{r}\deg(\lcm(\alpha_{i},\beta_{i+x-j}))\\&-\sum_{i=1}^{x-j}\deg(\beta_i)
  -\sum_{i=1}^{r}\max\{p_{i},q_{i+x-j}\}
  -\sum_{i=1}^{x-j}q_i\\
  = &\sum_{i=1}^{r}p_i+\sum_{i=1}^{x}c_i
  -\sum_{i=1}^{r}\max\{p_{i},q_{i+x-j}\}
  -\sum_{i=1}^{x-j}q_i.\\
\end{array}
$$
From (\ref{eqfromhifjioipol}) we obtain
$$
	\sum_{i=1}^{j}	\tilde a'_i\geq\sum_{i=1}^{x}c_i-\sum_{i=j+1}^{x}c_i=\sum_{i=1}^{j}c_i,  \quad 1\leq j \leq x-1. 
$$
Therefore,
$$
\min\{j\geq 1\; : \;\sum_{i=1}^{j}c_i>\sum_{i=i}^j \tilde a'_i\}
=x+1.
  $$
The result follows from 
 Theorem \ref{theoprescrifioi} (item \ref{ittheoprescrifioi2}).
 \hfill $\Box$
 
\medskip
The technique used in the following theorem is different from the one used in the previous results.
We first present a remark.

\begin{rem}\label{remsis}
Let $R(s)$ be a rational matrix with  $\tilde p_1, \dots, \tilde p_r$ as invariant orders at $\infty$. If the invariant rational functions of $R(\frac{1}{s})$ are
$\frac{\hat \eta_1(s)}{\hat \varphi_1(s)},\dots,   \frac{\hat \eta_r(s)}{\hat \varphi_r(s)}$, then
\begin{equation*}\frac{\hat \eta_i(s)}{\hat \varphi_i(s)}=s^{\tilde p_i}\frac{\hat \eta'_i(s)}{\hat \varphi'_i(s)}, \quad 1\leq i \leq r, \end{equation*}
where the triples of polynomials 
$(\hat  \eta'_i(s), \hat \varphi'_i(s), s)$
 are pairwise coprime for $1\leq i \leq r$ (see \cite[p. 724]{AnDoHoMa19} or \cite[Proposition 6.11]{AmMaZa15}).

\end{rem}

\begin{theorem} {\rm(Prescription of the infinite structure for  rational matrices)}
  \label{theoprescrioi}
 Let $R(s)\in\efe(s)^{m\times n}$ be a rational matrix, $\rank (R(s))=r$, and let 
 $\tilde p_1, \dots, \tilde p_r$ be its invariant orders at $\infty$.
 
	Let $z,x$ be integers such that $0\leq x\leq \min\{z, n-r\}$ and let $\tilde q_{1}\leq \dots \leq \tilde q_{r+x}$ be integers. There exists a rational matrix  $\widetilde W(s)\in \efe(s)^{z\times n}$ such that $\rank \left(\begin{bmatrix}R(s)\\\widetilde W(s)\end{bmatrix}\right)=r+x$ and
        $\begin{bmatrix}R(s)\\\widetilde W(s)\end{bmatrix}$ has 
 $\tilde q_1, \dots,  \tilde q_{r+x}$ as invariant orders at $\infty$
        if and only if (\ref{eqinterratioi}) holds.
        \end{theorem}
{\bf Proof.} The necessity follows from Theorem \ref{prescr4rat}.

For the sufficiency, assume that (\ref{eqinterratioi}) holds.
Let 
$\frac{\hat \eta_1(s)}{\hat \varphi_1(s)},\dots,   \frac{\hat \eta_r(s)}{\hat \varphi_r(s)}$ be the invariant rational functions of $R\left(\frac{1}{s}\right)$. We can write
$$\frac{\hat \eta_i(s)}{\hat \varphi_i(s)}=s^{\tilde p_i}\frac{\hat \eta'_i(s)}{\hat \varphi'_i(s)}, \quad 1\leq i \leq r,$$
where the triples of polynomials $ (\hat  \eta'_i(s), \hat \varphi'_i(s), s)$ are pairwise coprime (see Remark \ref{remsis}). Notice that 
$$
\begin{array}{lll}
  \hat \eta_{i}(s)=s^{\tilde p_i}\hat \eta'_i(s), & \hat \varphi_{i}(s)=\hat \varphi'_i(s), & \mbox{ if $\tilde p_i\geq 0$},\\
  \hat \eta_{i}(s)=\hat \eta'_i(s), & \hat \varphi_{i}(s)=s^{-\tilde p_i}\hat \varphi'_i(s), & \mbox{ if $\tilde p_i<0$},
\end{array}
$$
and
\begin{equation}\label{eqdivp}
\hat \eta'_i(s) \mid \hat \eta'_{i+1}(s), \quad
\hat \varphi'_{i+1}(s) \mid \hat \varphi'_{i}(s), \quad 1\leq i \leq r-1.
\end{equation}
For $1\leq i \leq x$, define
$$
\begin{array}{rlll}
  \hat \epsilon_i(s)=&s^{\tilde q_i}, &\hat \psi_i(s)=\hat \varphi'_1(s), & \mbox{ if $\tilde q_i\geq 0$},\\
  \hat \epsilon_i(s)=&1, & \hat \psi_i(s)=s^{-\tilde q_i}\hat \varphi'_1(s), & \mbox{ if $\tilde q_i<0$},
  \end{array}
  $$
  and for $1\leq i \leq r$,
$$
\begin{array}{rll}
  \hat \epsilon_{i+x}(s)=
s^{\tilde q_{i+x}}
  \hat \eta'_i(s), & \hat \psi_{i+x}(s)=\hat \varphi'_i(s),& \mbox{ if $\tilde q_{i+x}\geq 0$},\\
  \hat \epsilon_{i+x}(s)=\hat \eta'_i(s), & \hat \psi_{i+x}(s)=s^{-\tilde q_{i+x}}
  \hat \varphi'_i(s),&  \mbox{ if $\tilde q_{i+x}<0$}. 
  \end{array}
$$ Then,
$$
\begin{array}{ll}
  \frac{\hat \epsilon_i(s)}{\hat \psi_i(s)}=s^{\tilde q_i}\frac{1}{\hat \varphi'_1(s)}, & 1\leq i \leq x,
  \\
\frac{\hat \epsilon_{i+x}(s)}{\hat \psi_{i+x}(s)}=s^{\tilde q_{i+x}}\frac{\hat \eta'_i(s)}{\hat \varphi'_i(s)}, &1\leq i \leq r.
  \end{array}
$$
Since $\tilde q_i \leq \tilde q_{i+1}$,  $1\leq i \leq r+x-1$, from (\ref{eqdivp}) we obtain
$$
\hat \epsilon_i(s) \mid \hat \epsilon_{i+1}(s), \quad
\hat \psi_{i+1}(s) \mid \hat \psi_{i}(s), \quad 1\leq i \leq r+x-1,
$$
and from (\ref{eqinterratioi}) we obtain (recall that for $i>r+x$ we take $\epsilon_{i}(s)=0$ and $\psi_{i}(s)=1$)
 $$
\hat \epsilon_i(s)\mid \hat \eta_i(s)\mid \hat \epsilon_{i+z}(s),
\quad
\hat \psi_{i+z}(s)\mid \hat \varphi_i(s)\mid \hat \psi_{i}(s), \quad 1\leq i\leq r.
$$
By Theorem   \ref{lemmainvrat}, there is a  rational matrix  $\widehat W(s)\in \efe(s)^{z\times n}$,  $\rank \left(\begin{bmatrix}R(\frac{1}{s})\\\widehat  W(s)\end{bmatrix}\right)=r+x$, such that
$\begin{bmatrix}R(\frac{1}{s})\\\widehat  W(s)\end{bmatrix}$ has $\frac{\hat \epsilon_1(s)}{\hat \psi_1(s)}, \dots, \frac{\hat \epsilon_{r+x}(s)}{\hat \psi_{r+x}(s)}$ as invariant rational functions.

Let $\widetilde W(s)=\widehat W(\frac{1}{s})$ and $\widetilde Q(s)=\begin{bmatrix}R(s)\\\widetilde  W(s)\end{bmatrix}$. Then,   $\widetilde Q(\frac{1}{s})=\begin{bmatrix}R(\frac{1}{s})\\\widehat  W(s)\end{bmatrix}$  and
$\tilde q_1, \dots, \tilde q_{r+x}$ are the invariant orders at $\infty$ of $\widetilde Q(s)$ (Remark \ref{remsis}).
\hfill $\Box$

\medskip
The next example shows the difference between the rational and polynomial cases when prescribing the infinite structure.

\begin{example}\label{exprioi}
  Let
  $P(s)=\begin{bmatrix}s&0\end{bmatrix}\in \FF[s]^{1\times 2}$. The matrix $P(s)$ has $\alpha_1(s)=s$ as invariant factor, $p_1=-1$ as invariant order at $\infty$ and $c_1=0$ as column minimal index.

  Let $z=x=1$, and $q_1=-1, q_2=+1$. Then  (\ref{eqinterpolioi}) holds, but
  (\ref{eqfromhifjioipol}) is not satisfied. Therefore, there is no polynomial matrix $Q(s)=\begin{bmatrix}P(s)\\W(s)\end{bmatrix}\in \FF[s]^{2\times 2}$ of $\rank Q(s)=2$  with  $q_1=-1, q_2=+1$ as invariant orders at $\infty$.
   If there were such  a polynomial matrix,  then by Theorem \ref{theoprescr4ioi}, the invariant factors $\beta_1(s), \beta_2(s)$  of $Q(s)$ would  satisfy $\alpha_1(s)=s\mid \beta_2(s)$ and  $\deg(\beta_1)+\deg( \beta_2)=0$, which leads to a contradiction.

  However, if we allow the completion to be rational, it is possible to obtain the desired invariants.
For example, 
the rational matrix   $\widetilde Q(s)=\begin{bmatrix}
s & 0\\
0 & \frac{1}{s}\end{bmatrix}\in \FF(s)^{2\times 2}$
has  $\tilde q_1=-1, \tilde q_2=+1$ as invariant orders at $\infty$.
\end{example}

As mentioned, a solution to the row completion problem for polynomial   matrices when the finite structure is prescribed follows from  Theorem \ref{lemmaSaTh}.
In this theorem no condition is imposed on the invariant orders at $\infty$, and therefore on the degree of the completed matrix $Q(s)$. To prescribe the degree of $Q(s)$, we must prescribe the first order at $\infty$ of $Q(s)$, $q_1=-\deg(Q(s))$, which shall  satisfy $q_1=-\deg(Q(s))\leq-\deg(P(s))=p_1$.

Theorem \ref{lemmaSaTh} was later generalized to rational matrices in Theorem 
\ref{lemmainvrat}. In the next theorem, we give a solution to the  row completion problem for rational   matrices when the finite structure and the first invariant order at $\infty$ of the completed matrix is prescribed.

\begin{theorem} {\rm (Prescription of the finite structure and the first invariant order at $\infty$)}
  \label{theoprescrifdegree}
 Let $R(s)\in\efe(s)^{m\times n}$ be a rational matrix, $\rank (R(s))=r$. 
 Let $\frac{\eta_1(s)}{\varphi_1(s)},\dots, \frac{\eta_r(s)}{\varphi_r(s)}$ 
	be its invariant rational functions,  
$\tilde p_1, \dots, \tilde p_r$ its invariant orders at $\infty$, 
	$\bc=(c_1,  \dots,  c_{n-r})$  its
	column minimal indices.
	
	Let $z,x$ be integers such that $0\leq x\leq \min\{z, n-r\}$, let $\epsilon_1(s)\mid \dots \mid\epsilon_{r+x}(s)$ and 
        $\psi_{r+x}(s)\mid \dots \mid\psi_{1}(s)$ be monic polynomials such that $\frac{\epsilon_i(s)}{\psi_i(s)}$ are irreducible rational functions, $1\leq i \leq r+x$, and let 
        $\tilde q_{1}$ be an integer
         $\tilde q_{1}\leq\tilde p_{1}$. 
       There exists a
 rational matrix  $\widetilde W(s)\in \efe(s)^{z\times n}$ such that $\rank \left(\begin{bmatrix}R(s)\\\widetilde W(s)\end{bmatrix}\right)=r+x$ and
        $\begin{bmatrix}R(s)\\\widetilde W(s)\end{bmatrix}$ has $\frac{\epsilon_1(s)}{\psi_1(s)}, \dots, \frac{\epsilon_{r+x}(s)}{\psi_{r+x}(s)}$ as invariant rational functions and 
 $\tilde q_1$ as first invariant order at $\infty$
        if and only if
        (\ref{eqinterratnum}),
         (\ref{eqinterratden}) and 
\begin{equation}\label{eqfromhifjioiq1rat}
\begin{array}{ll}
            \sum_{i=1}^{r}\Delta\left(\frac{\eta_i}{\varphi_i}, \frac{\epsilon_{i+x-j}}{\psi_{i+x-j}}\right)
            +\sum_{i=1}^{x-j}\Delta\left(\frac{\epsilon_{i}}{\psi_{i}}\right)-
  \sum_{i=1}^{r}\Delta\left(\frac{\eta_{i}}{\varphi_{i}}\right) \\\leq          
            \sum_{i=j+1}^{x}c_i+(j-x)\tilde q_1,\quad
  0\leq j \leq x-1.
  \end{array}
  \end{equation}
  \end{theorem}

{\bf Proof.} 
The proof is similar to that of Theorem
\ref{theoprescrioipol} and 
\cite[Theorem 4.12]{AAmBaMaRo24}.  We also precise the steps.

Let $\bu=(u_1, \dots,  u_{m-r})$ be the row minimal indices of $R(s)$.
Assume that 
there exists  $\widetilde W(s)\in \efe(s)^{z\times n}$ such that $\rank\left(\begin{bmatrix}R(s)\\\widetilde W(s)\end{bmatrix}\right)=r+x$ and $\begin{bmatrix}R(s)\\\widetilde W(s)\end{bmatrix}$ has
 $\frac{\epsilon_1(s)}{\psi_1(s)}, \dots, \frac{\epsilon_{r+x}(s)}{\psi_{r+x}(s)}$ as invariant rational functions  and 
  $\tilde q_1\le \dots \le \tilde q_{r+x}$ as invariant orders at  $\infty$.

By Theorem \ref{theoprescrifioi} and Remark \ref{4.1819-4.17} we obtain (\ref{eqinterratnum})--(\ref{eqinterratioi}) and 
        \begin{equation}\label{eqfromhifjioipol1deg}
          \begin{array}{l}
          \sum_{i=1}^{r}\Delta\left(\frac{\eta_i}{\varphi_i}, \frac{\epsilon_{i+x-j}}{\psi_{i+x-j}}\right)
            +\sum_{i=1}^{x-j}\Delta\left(\frac{\epsilon_{i}}{\psi_{i}}\right)
+\sum_{i=1}^{r}\max\{ \tilde p_i,  \tilde q_{i+x-j}\}\\+\sum_{i=1}^{x-j} \tilde q_i
  +\sum_{i=1}^{m-r}u_i+\sum_{i=1}^{j}c_i+\sum_{i=x+1}^{n-r}c_i\leq  0,
  \quad  0\leq j \leq x-1.
  \end{array}
        \end{equation}
We have 
$$\sum_{i=1}^{r}\max\{ \tilde p_i,  \tilde q_{i+x-j}\}\\+\sum_{i=1}^{x-j} \tilde q_i\geq 
\sum_{i=1}^{r} \tilde p_i+(x-j) \tilde q_1, \quad 0\leq j \leq x-1.$$
 Thus, from  (\ref{eqfromhifjioipol1deg}) we obtain
$$\begin{array}{l}
      \sum_{i=1}^{r}\Delta\left(\frac{\eta_i}{\varphi_i}, \frac{\epsilon_{i+x-j}}{\psi_{i+x-j}}\right)
            +\sum_{i=1}^{x-j}\Delta\left(\frac{\epsilon_{i}}{\psi_{i}}\right)
+\sum_{i=1}^{r} \tilde p_i\\
  +\sum_{i=1}^{m-r}u_i+\sum_{i=1}^{j}c_i+\sum_{i=x+1}^{n-r}c_i\leq  (j-x)\tilde q_1,       
 \quad  0\leq j \leq x-1,
  \end{array}
  $$
  which by Theorem \ref{theoexistencerat} is equivalent to (\ref{eqfromhifjioiq1rat}).

Conversely, 
assume that 
(\ref{eqinterratnum}),
         (\ref{eqinterratden}) and 
(\ref{eqfromhifjioiq1rat})  hold.
Let $$
\tilde t=
\sum_{i=1}^{r}\Delta\left(\frac{\eta_i}{\varphi_i}\right)- 
         \sum_{i=1}^{r}\Delta\left(\frac{\eta_i}{\varphi_i},
         \frac{\epsilon_{i+x}}{\psi_{i+x}}
         \right)-   
\sum_{i=1}^{x}\Delta\left(
         \frac{\epsilon_{i}}{\psi_{i}}
         \right)+ \sum_{i=1}^x c_i-x\tilde q_1.
$$
If $x=0$, from 
(\ref{eqinterratnum}) and 
(\ref{eqinterratden}) 
we obtain $\tilde t=0$. If $x>0$, from (\ref{eqfromhifjioiq1rat})
we have $\tilde t\geq 0$.

Define
$$ 
\begin{array}{lcll}
\tilde q_i&=&\tilde q_1, &  1\leq i \leq x,\\
  \tilde q_{i+x}&=&\tilde p_{i}, & 1\leq i \leq r-1,\\
  \tilde q_{r+x}&=&\tilde p_{r}+ \tilde t.
  \end{array}
  $$
  As $\tilde q_1 \leq \tilde p_1$ and $\tilde t\geq 0$, we have $\tilde q_1\leq \dots \leq\tilde q_{r+x}$, and (\ref{eqinterratioi}) holds.

If $x>1$, 
for $1\leq j \leq x$, we have $\tilde q_{i+x-j} \leq \tilde q_{i+x-1}\leq \tilde p_i$, $1\leq i \leq r$, therefore 
$$
\sum_{i=1}^{r}\max\{\tilde p_i, \tilde q_{i+x-j}\}+\sum_{i=1}^{x-j}\tilde q_i=\left\{\begin{array}{l}
\sum_{i=1}^{r}\tilde p_i+\tilde t+ x\tilde q_1, \quad  j=0,\\
\sum_{i=1}^{r}\tilde p_i+(x-j) \tilde q_1, \quad 1\leq j  \leq x.
\end{array}\right. 
$$
Thus, from Theorem
\ref{theoexistencerat} 
 we obtain
\begin{equation}\label{eqcero2}
 \begin{array}{rl}
          &\sum_{i=1}^{r}\Delta\left(\frac{\eta_i}{\varphi_i}, \frac{\epsilon_{i+x}}{\psi_{i+x}},\tilde p_i, \tilde q_{i+x}\right)
            +\sum_{i=1}^{x}\Delta\left(\frac{\epsilon_{i}}{\psi_{i}}, \tilde q_i\right)
  +\sum_{i=1}^{m-r}u_i+\sum_{i=x+1}^{n-r}c_i\\
  = &\sum_{i=1}^{r}\Delta\left(\frac{\eta_i}{\varphi_i}, \frac{\epsilon_{i+x}}{\psi_{i+x}}\right)
            +\sum_{i=1}^{x}\Delta\left(\frac{\epsilon_{i}}{\psi_{i}}\right)+
            \tilde t+ x\tilde q_1-\sum_{i=1}^r \Delta\left(\frac{\eta_i}{\varphi_i}\right)-\sum_{i=1}^{x}c_i\\=&0,
  \end{array}
  \end{equation}
  and from Theorem
\ref{theoexistencerat} and (\ref{eqfromhifjioiq1rat}) 
$$\begin{array}{l}
  \sum_{i=1}^{r}
  \Delta\left(\frac{\eta_i}{\varphi_i}, \frac{\epsilon_{i+x-j}}{\psi_{i+x-j}},\tilde p_i, \tilde q_{i+x-j}\right)
            +\sum_{i=1}^{x-j}\Delta\left(\frac{\epsilon_{i}}{\psi_{i}}, \tilde q_i\right)
  +\sum_{i=1}^{m-r}u_i+\sum_{i=1}^{j}c_i\\+\sum_{i=x+1}^{n-r}c_i=
  \sum_{i=1}^{r}
  \Delta\left(\frac{\eta_i}{\varphi_i}, \frac{\epsilon_{i+x-j}}{\psi_{i+x-j}}\right)
            +\sum_{i=1}^{x-j}\Delta\left(\frac{\epsilon_{i}}{\psi_{i}}\right)+
            \sum_{i=1}^{r}\tilde p_i+(x-j) \tilde q_1\\
  +\sum_{i=1}^{m-r}u_i+\sum_{i=1}^{j}c_i+\sum_{i=x+1}^{n-r}c_i=
  \sum_{i=1}^{r}
  \Delta\left(\frac{\eta_i}{\varphi_i}, \frac{\epsilon_{i+x-j}}{\psi_{i+x-j}}\right)
            +\sum_{i=1}^{x-j}\Delta\left(\frac{\epsilon_{i}}{\psi_{i}}\right)\\-
\sum_{i=1}^r \Delta\left(\frac{\eta_i}{\varphi_i}\right)-\sum_{i=j+1}^{x}c_i+           
            (x-j) \tilde q_1\leq 0, \quad 1\leq j \leq x-1.
\end{array}
$$

If $x<z$ or $x=z=n-r$ the result follows from  Theorem \ref{theoprescrifioi} (item \ref{ittheoprescrifio11}).

If $x=z<n-r$, 
let
$\tilde \ba'=(\tilde a'_1, \dots, \tilde a'_x)$
be defined as in (\ref{eqdeftildeap}).
From (\ref{eqinterratnum}), (\ref{eqinterratden}), (\ref{eqinterratioi}) and (\ref{eqcero2}) we have
$$\begin{array}{rl}
\sum_{i=1}^{x}	\tilde a'_i=&
\sum_{i=1}^{r+x}
\Delta\left(\frac{\epsilon_{i}}{\psi_{i}}, \tilde q_{i}\right)
            -\sum_{i=1}^{r}\Delta\left(\frac{\eta_{i}}{\varphi_{i}}, \tilde p_i\right)\\&=
  \tilde t + x\tilde q_1+\sum_{i=1}^{r+x}\Delta\left(\frac{\epsilon_{i}}{\psi_{i}}\right)
            -\sum_{i=1}^{r}\Delta\left(\frac{\eta_{i}}{\varphi_{i}}\right)\\&=
            \sum_{i=1}^{x}c_i-\sum_{i=1}^{r}\Delta\left(\frac{\eta_i}{\varphi_i}, \frac{\epsilon_{i+x}}{\psi_{i+x}}\right)+\sum_{i=1}^{r}\Delta\left(\frac{\epsilon_{i+x}}{\psi_{i+x}}\right)=
             \sum_{i=1}^{x}c_i.
\end{array}
$$

Let $j\in\{1,\dots, x-1\}$. Then 
$$
\begin{array}{rl}
	\sum_{i=1}^{j}	\tilde a'_i=&
  \sum_{i=1}^{r+x}
  \Delta\left(\frac{\epsilon_{i}}{\psi_{i}}, \tilde q_{i}\right)-
\sum_{i=1}^{r}
  \Delta\left(\frac{\eta_i}{\varphi_i}, \frac{\epsilon_{i+x-j}}{\psi_{i+x-j}},\tilde p_i, \tilde q_{i+x-j}\right)
            \\&-\sum_{i=1}^{x-j}\Delta\left(\frac{\epsilon_{i}}{\psi_{i}}, \tilde q_i\right)=
            \sum_{i=1}^{r}
  \Delta\left(\frac{\eta_i}{\varphi_i}, \tilde p_i\right)+\sum_{i=1}^{x}c_i\\&
-
\sum_{i=1}^{r}
  \Delta\left(\frac{\eta_i}{\varphi_i}, \frac{\epsilon_{i+x-j}}{\psi_{i+x-j}},\tilde p_i, \tilde q_{i+x-j}\right) -\sum_{i=1}^{x-j}\Delta\left(\frac{\epsilon_{i}}{\psi_{i}}, \tilde q_i\right)\\
  =&
  \sum_{i=1}^{r}
  \Delta\left(\frac{\eta_i}{\varphi_i}\right)+\sum_{i=1}^{x}c_i
-
\sum_{i=1}^{r}
  \Delta\left(\frac{\eta_i}{\varphi_i}, \frac{\epsilon_{i+x-j}}{\psi_{i+x-j}}\right)\\& -\sum_{i=1}^{x-j}\Delta\left(\frac{\epsilon_{i}}{\psi_{i}}\right)+(j-x)\tilde q_1.\\
\end{array}
$$
From (\ref{eqfromhifjioiq1rat}) we obtain
$$
	\sum_{i=1}^{j}	\tilde a'_i\geq\sum_{i=1}^{x}c_i-\sum_{i=j+1}^{x}c_i=\sum_{i=1}^{j}c_i,  \quad 1\leq j \leq x-1. 
$$
Therefore,
$$
\min\{j\geq 1\; : \;\sum_{i=1}^{j}c_i>\sum_{i=i}^ja'_i\}
=x+1.
  $$
The result follows from 
 Theorem \ref{theoprescrifioi} (item \ref{ittheoprescrifioi2}).
 \hfill $\Box$

\bibliographystyle{acm}
\bibliography{references}

\end{document}